\documentclass[12pt]{article}

    \usepackage{amsmath}
    \usepackage{amsfonts,color,amssymb,mathrsfs,stmaryrd}

\usepackage{graphicx}
\usepackage{wrapfig}
\usepackage{mathrsfs}

    \newcommand{\eq}[1]{$(\ref{#1})$}

    \newcommand{\eps}{\varepsilon}

\newcommand{\PP}{\mathbb{P}}
\newcommand{\RR}{\mathbb{R}}

\renewcommand{\H}{{\cal P}}

\newcommand{\X}{{\cal X}}
\newcommand{\0}{{\bf 0}}

    \newtheorem{theo}{Theorem}
    \newtheorem{coro}{Corollary}

    \newtheorem{lemm}{Lemma}
    \newtheorem{defn}{Definition}

        \def\N{\mathbb{N}}

    \def\E{\mathbb{E}}
    
    \def\0{{\bf 0}}
    
    \def\R{\mathbb{R}}
    \def\RR{\mathbb{R}}
    \def\PP{\mathbb{P}}

    \def\Ed{{\cal E}}

    \def\G{{\cal G}}

    \renewcommand{\E}{\mathbb E \,}
    \newcommand{\C}{{\cal C}}

    \newcommand{\tod}{\stackrel{{\cal D}}{\longrightarrow}}

     \newcommand{\liml}{\lim_{\la \to \infty} }

    \newcommand{\eqco}{\setcounter{equation}{0}}
    \newcommand{\thco}{\setcounter{theo}{0}}
    \newcommand{\prco}{\setcounter{prop}{0}}
    \newcommand{\laco}{\setcounter{lemm}{0}}
    \newcommand{\coco}{\setcounter{coro}{0}}
    \newcommand{\cjco}{\setcounter{conj}{0}}
    
    \newcommand{\deco}{\setcounter{defn}{0}}
    \newcommand{\allco}{\eqco  \thco \prco \laco \coco \cjco \deco}
\newcommand{\M}{{\cal M}}

    \renewcommand{\P}{{{\cal P}}}

    \newcommand{\A}{{\cal A}}
    \newcommand{\Cov}{{\rm Cov}}

    \newcommand{\Var}{{\rm Var}}
    \newcommand{\var}{{\rm Var}}

    \newcommand{\card}{{\rm card}}
    \newcommand{\Vol}{{\rm Vol}}
    \newcommand{\diam}{{\rm diam}}

        \newcommand{\F}{{\cal F}}

    \def\bdm{\begin{displaymath}}
    \newcommand{\edm}{\end{displaymath}}
    \def\benu{\begin{enumerate}}
    \def\eenu{\end{enumerate}}
    \def\beqn{\begin{equation}}
    \def\eeqn{\end{equation}}
    \def\be{\begin{equation}}
    \def\ee{\end{equation}}
    \def\bea{\begin{eqnarray}}
    \def\eea{\end{eqnarray}}
    \newcommand{\bean}{\begin{eqnarray*}}
    \newcommand{\eean}{\end{eqnarray*}}
    \newcommand{\bear}{\begin{eqnarray}}
    \newcommand{\eear}{\end{eqnarray}}

    \def\R{\mathbb{R}}
    \def\D{\mathbb{D}}

    \def\A{{\cal A}}
    \def\la{{\lambda}}
    \def\ka{{\kappa}}
    \def\b{{\beta}}

    \def\qed{\hfill\hbox{${\vcenter{\vbox{
        \hrule height 0.4pt\hbox{\vrule width 0.4pt height 6pt
        \kern5pt\vrule width 0.4pt}\hrule height 0.4pt}}}$}}

\def\setminusx{{\setminus\{x\}}}
\def\setminusz{{\setminus\{z\}}}

\def\bs{{\bf s}}

\def\ignore#1{}
\def\Ref#1{(\ref{#1})}
\def\tworho{\rho} 
\def\fourrho{2\rho}
\def\eightrho{4\rho}

\begin{document}
\title{\bf Normal approximation of Gibbsian sums in geometric probability}

    \author{Aihua Xia,   J. E. Yukich}

    \date{\today}
    \maketitle

 \footnotetext{Research supported in part by ARC Discovery Grant DP130101123 (AX) and NSF grants
    DMS-1106619, DMS-1406410 (JY)}

 \begin{abstract}
This paper concerns the asymptotic behavior of a random variable $W_\la$ resulting from
the summation of the functionals of a Gibbsian spatial point process over windows $Q_\la \uparrow \R^d$. 
 We establish conditions ensuring that $W_\la$ has volume order fluctuations, that is they coincide with the fluctuations
 of functionals of Poisson spatial point processes.  We combine this result with Stein's method to deduce 
rates of normal approximation for $W_\la$, as $\la\to\infty$.
Our general results  
establish variance asymptotics and central limit theorems for statistics of random geometric and
related Euclidean graphs on Gibbsian input.  We also establish similar limit theory for
 claim sizes of insurance models with Gibbsian input, the number of maximal points of a Gibbsian sample,
and the size of spatial birth-growth models with Gibbsian input.



\end{abstract}

\section{Introduction and main results}
 \allco
Functionals of large  geometric structures on finite input $\X \subset \R^d$ often consist
of sums of spatially dependent terms admitting the representation
\be \label{basic}\sum_{x \in \X} \xi(x,\X), \ee  where the $\R^+$-valued {\em score function} $\xi$, defined on pairs $(x, \X)$,   represents the {\em interaction}
of $x$ with respect to $\X$.  The sums \eqref{basic} typically describe a global feature of an underlying  geometric property in terms of a sum of local contributions $\xi(x,\X)$.

A large and diverse number of functionals and statistics
in stochastic geometry, applied geometric probability, and spatial statistics may be
cast in the form \eqref{basic} for appropriately chosen $\xi$.  The behavior
of these statistics  on random input $\X$  can  be deduced from general limit theorems \cite{BY, PeEJP, PeBer, PY4, PY5} for \eqref{basic}
provided $\X$ is either a Poisson or binomial point process.  This has led to solutions of
problems in random sequential packing \cite{PY2}, random graphs \cite{PeEJP, PeBer, PY1, PY4, Wa}, percolation models \cite{LP},
analysis of data on manifolds \cite{PY6},  and
convex hulls of i.i.d. samples \cite{CSY, CY, CY2}, among others.

When $\X$ is neither Poisson nor binomial input, the limit theory of  \eqref{basic} is less well understood.
Our main purpose is to redress this for Gibbsian input.
For all $\la \in [1, \infty)$ consider the functionals
$$W_\la:= \sum_{x \in {\H}_\la^{\b \Psi}} \xi(x,  {\H}_\la^{\b \Psi} \setminus \{x\}),$$
where
${\H}_\la^{\b \Psi}$
is 
the restriction of a Gibbs point process  ${\H}^{\b \Psi}$ on $\R^d$ to
$Q_\la:= [-\la^{1/d}/2, \la^{1/d}/2]^d$.  The process ${\H}^{\b \Psi}$
has potential
$\Psi$, it is absolutely continuous with respect to
 a reference homogeneous Poisson point process $\tilde\P_{\tau}$ of intensity $\tau$,
and $\b$ is the inverse temperature.
In general, even for the simplest of score functions $\xi$, the Gibbsian functional $W_\la$  may neither  enjoy asymptotic normality nor have volume order fluctuations, i.e.,
 { $\Var W_\la$  may not be of order} $\Vol (Q_\la)$; see \cite{MY}.  On the other hand, if both the Gibbsian input and the score function have rapidly decaying spatial
dependencies, then one could expect that $W_\la$ behaves like a sum of i.i.d. random variables.

 We have three goals.  The first is to show that given a potential $\Psi$,
 there is a range of inverse temperature and intensity parameters
 $\b$ and $\tau$ such that for any locally determined score function, the Gibbsian functional $ W_\la$ has volume order fluctuations.
In other words,  the fluctuations for $ W_\la$ coincide with
 those when ${\H}_\la^{\b \Psi}$ is replaced by Poisson or binomial input.
 This strengthens the central limit theorems of \cite{SY}, which
 depend crucially on volume order fluctuations.
 Our second goal is to prove a rate of convergence to the normal for $ (W_\la - \E W_\la)/ \sqrt{ \Var W_\la }$ for general score functions $\xi$, including
 those which are non-translation invariant.
Formal statements of these  results are given in Theorems \ref{main0}-\ref{main-noti}.
Thirdly, we use our general results to deduce rates of normal convergence for
(i) statistics of random geometric and  Euclidean graphs on Gibbsian input, (ii)  the number of
claims in an insurance model with claim locations and times given by Gibbsian input, (iii) the
number of maximal points in a Gibbs sample,  as well as (iv) functionals of spatial birth-growth models with Gibbsian input.
This extends the  central limit theorems and second order results of  \cite{BX1, LPS, PeEJP, PY4, PY5} to Gibbsian input.


\vskip.5cm

\subsection{Notation and terminology}

\noindent (i)\ {\bf Gibbs point processes.}
Quantifying spatial dependencies of Gibbs point processes is difficult in general. However spatial dependencies readily become transparent
 when a Gibbs point process is viewed as an algorithmic construct. 
 As shown in \cite{SY}, this is feasible whenever
$\Psi$ belongs  to the class of potentials $\bf{\Psi}^*$ containing pair potentials, continuum
 Widom-Rowlinson potentials,  area interaction potentials,  hard core potentials
 and potentials generating a truncated Poisson point process.

We review the algorithmic construction of
Gibbs point processes developed  in \cite{SY}, and inspired by \cite{FFG}. Define for $\Psi \in \mathbf{\Psi^*}$ and finite $\X \subset \R^d$ the local energy function
 $$
 \Delta^\Psi(\0, \X):= \Psi(\X \cup \{\0\} ) - \Psi(\X), \ \0 \notin \X.
 $$
 Here $\0$ denotes a point at the origin of $\R^d$.
 Proposition 2.1~(i) of \cite{SY} shows that for $\X \subset \R^d$ locally finite,
 \be \label{Dell}
 \Delta^\Psi(\0, \X):= \lim_{r \to \infty} \Delta^\Psi(\0, \X \cap B_r(\0))
 \ee
 is well-defined, where $B_r(x):=\{y:\ |x - y| \le r\}$ is the Euclidean ball with center $x$ and radius $r$.
 $\Psi$ has {\em finite or bounded range}
 if there is $r^{\Psi} \in (0, \infty)$ such that for all
 finite $\X \subset \R^d$ we have $\Delta^\Psi(\0, \X) = \Delta^\Psi(\0, \X \cap B_{r^{\Psi}}(\0)  ).$
 {With the exception of the pair potential, all potentials in $\mathbf{\Psi^*}$ have finite range (Lemma 3.1 of \cite{SY}).}
 For such $\Psi$ we put
 $$
 m_0^{\Psi}:= \inf_{ \X \ {\rm{locally\ finite}} } \Delta^\Psi(\0, \X)
 $$
 and
 \be \label{defRP}
 {\cal R}^\Psi:= \{(u,v) \in (\R^+)^2: \ uv_d \exp(-v m_0^{\Psi}) (r^{\Psi} +1)^d < 1 \},
 \ee
where $v_d:= \pi^{d/2}[ \Gamma(1 + d/2)]^{-1}$ is
the volume of the unit ball in $\R^d.$  When $\Psi$ is a pair potential,  then the factor
 $(r^{\Psi} +1)^d$ in \eqref{defRP} is replaced by the moment of
 an exponentially decaying random variable as in  (3.7)
 of \cite{SY}.

Let $(\varrho(t))_{t \in \R}$ be a stationary homogeneous free
birth and death process on $\R^d$ with these dynamics:
\begin{itemize}
  \item A new point $x \in \R^d$ is born in $\varrho_t$ during the time interval $[t-dt,t]$
        with probability $\tau dx dt,$
  \item An existing point $x \in \varrho_t$ dies during the time interval
        $[t-dt,t]$ with probability $dt,$ that is
        the lifetimes of points of the process are independent standard
        exponential.
 \end{itemize}
The unique stationary and reversible measure for this process is
the law of the Poisson point process $\tilde{\cal P}_\tau.$


Following \cite{SY}, for each $\Psi \in \mathbf{\Psi^*},$ we use a dependent thinning procedure on
$(\varrho(t))_{t \in \R}$ to algorithmically construct a Gibbs point process  $\P^{\b \Psi}$ on $\R^d$, one whose
law is absolutely continuous  with respect to the reference  point process $\tilde{\cal P}_\tau.$ Section 3 recalls  some of the salient properties
of $\P^{\b \Psi}$ .

For arbitrary $(\tau,\beta)$ and arbitrary $\Psi$, the asymptotic behavior of  $W_\la$ may involve non-standard scaling and non-standard limits.
However, if $\P^{\b \Psi}$ is {\em admissible in the sense that $(\tau, \beta) \in {\cal R}^\Psi$ and $\Psi \in \mathbf{\Psi^*},$
then we shall show that $W_\la$  behaves like a classical sum of i.i.d. random variables.} Henceforth, and
without further mention,  we shall always assume that $\P^{\b \Psi}$ is admissible.
Recall that $Q_\la:= [-\la^{1/d}/2, \la^{1/d}/2]^d$ and put $Q_\infty := \R^d$. Given  $\la \in [1, \infty]$, $\Psi \in \mathbf{\Psi^*},$ and $(\tau, \beta) \in {\cal R}^\Psi$,
we let
  \be \label{Gibbs} {\H}_\la^{\b \Psi} := \P^{\b \Psi} \cap Q_\la. \ee
By convention we have  ${\H}_{\infty}^{\b \Psi}:=  {\H}^{\b \Psi}.$

 \vskip.3cm
\noindent {(ii) \bf Poisson-like point processes.} 
 A point process $\Xi$ on
${\mathbb R}^d$ is {\it stochastically dominated} by the reference process
$\tilde{\cal P}_\tau$ if for all Borel sets $B \subset \R^d$ and  $n \in \N$ we
have $\mathbb \PP[\text{card}(\Xi \cap B) \geq n] \leq \mathbb
\PP[\text{card}(\tilde{\cal P}_\tau \cap B) \geq n]$. 
As in \cite{SY}, we  say that $\Xi$ is {\em Poisson-like} if (a) $\Xi$ is stochastically
dominated by $\tilde\P_\tau$ and (b) there exists ${c} \in (0, \infty)$ and $r_1 \in (0, \infty)$ such that for all $r \in (r_1, \infty)$, $x \in \R^d$, and any point set ${\cal E}_r(x)$ in $B_r^c(x)$,
the conditional probability
of  $B_r(x)$ not being hit by $\Xi$, given that $\Xi \cap B_r(x)^c$
coincides with ${\cal E}_r(x)$,
satisfies
 \be \label{LBd}
    {\mathbb P}[\Xi \cap B_r(x) = \emptyset \ | \{  \Xi \cap B_r(x)^c = {\cal E}_r(x)  \} ] \leq \exp(-{c} r^d). 
 \ee
Poisson-like processes have void probabilities analogous to those of
homogeneous Poisson processes, justifying the choice of terminology.
Lemma 3.3 of \cite{SY} shows that admissible Gibbs  processes
$\P^{\b \Psi}$ are Poisson-like.

\ \

\noindent{(iii) \bf Translation invariance}.  $\xi$ is {\em translation invariant} if for all $x \in \R^d$ and locally finite $\X \subset \R^d$ we have $\xi(x, \X) = \xi(x + y, \X + y)$ for all $y \in \R^d$.

\ \

\noindent{(iv) \bf Moment conditions.}  Let $\|X\|_q$ denote the $q$ norm of the random variable $X$. Say that $\xi$ satisfies the
$q$-moment condition if
\be \label{mom} w_q:= \sup_{\la \in [1, \infty] } \sup_{x \in Q_\la} \| \xi(x,{\H}_\la^{\b \Psi} \cup \{x\} )\|_q < \infty.
\ee


\noindent{(v) \bf Stabilization.}  Given a locally
finite point set $\X$, write $\X^z$ for $\X \cup
\{z\}$ if $z \in \R^d$ and $\X^z=\X$ if $z=\emptyset$.  The following definition of stabilization is similar to that in \cite{BX, PeEJP, PeBer, PY4, PY5}
except now we consider Gibbsian input, instead of Poisson or binomial input.
\begin{defn} \label{def5}   $\xi$ is a
 {\it stabilizing functional} with respect to the  Poisson-like process
 $\Xi$ if for all $x \in \R^d$, all $z \in \R^d \cup \{\emptyset\}$,  and almost all realizations $\X$ of $\Xi$
  there exists $R :=  R^{\xi}(x,\X^z) \in (0, \infty)$ (a `radius of stabilization')  such that
 \begin{equation}\label{StabRelW}
  \xi(x,\X^z  \cap B_R(x)) = \xi(x, (\X^z  \cap B_R(x)) \cup {\cal Y})
 \end{equation}
 for all locally finite point sets ${\cal Y} \subseteq \R^d   \setminus B_R(x).$
\end{defn}

Stabilization  of $\xi$ on
$\Xi$ implies that $\xi(x, \X^z)$ is wholly determined by the point
configuration $\X^z \cap B_{R^\xi}(x)$. It also yields $\xi(x,\X^z
\cap B_r(x)) = \xi(x,\X^z \cap B_{R^\xi}(x))$ for $r \in [R^\xi,
\infty)$. Stabilizing functionals  can thus be a.s.
extended to the entire process
 $\Xi^z$, that is to say for all $x \in \R^d$ and $z \in \R^d \cup \{\emptyset\}$
 we have
\begin{equation} \label{2.2a}
 \xi(x,\Xi^z) := \lim_{r\to\infty} \xi(x,\Xi^z \cap B_r(x)) \  \ \
\text{a.s.}  \end{equation}

Given $s > 0$ and any simple point process $\Xi$, including Poisson-like processes, define the conditional tail probability
$$
t(\Xi,s):=  \sup_{x  \in \R^d} \sup_{z \in \R^d \cup
\{\emptyset\}} \PP [R^{\xi}(x, \Xi^z) > s|\Xi\{x\}=1].$$
The conditional distribution of $\Xi$ given that $\Xi\{x\}=1$ is the Palm distribution of $\Xi$ at $x$ \cite[Chapter~10]{Kallenberg83} and the conditional probability can be intuitively interpreted as
$$\PP [R^{\xi}(x, \Xi^z) > s|\Xi\{x\}=1]=\lim_{\epsilon\downarrow0}\PP [\sup_{y\in B_\epsilon(x)\cap \Xi} R^{\xi}(y, \Xi^z) > s|\Xi(B_\epsilon(x))=1].
$$
We say that  $\xi$ is {\it stabilizing in the wide sense}
 if for every Poisson-like  process $\Xi$ we have $t(\Xi,s) \to 0$ as $s \to \infty.$
 Further, $\xi$ is {\it exponentially
stabilizing in the wide sense}
 if for every Poisson-like  process $\Xi$
 we have
 \be \label{ExpStab} \limsup_{s \to \infty} s^{-1} \ln t(\Xi,s) < 0. \ee
Exponential stabilization of $\xi$ with respect to the augmented point set  $\Xi^z$ ensures that covariances of scores at points $x$ and $y$, as given at \eqref{defc2x}, decays exponentially fast with $|x - y|$, implying that
$W_\la$ has at most volume order fluctuations, as seen in the proof of Lemma~\ref{varld}.
\noindent
Notice that for $\la$ large we  have $ R^{\xi}(x, \Xi^z \cap Q_\la) \leq R^{\xi}(x, \Xi^z)$
and thus \eqref{ExpStab}  holds with   $t(\Xi,s)$ replaced by
\be \label{ExpStabA}
\limsup_{\la \to \infty} \sup_{x  \in Q_\la} \sup_{z \in \R^d \cup
\{\emptyset\}} \PP [R^{\xi}(x, \Xi^z \cap Q_\la) > s|\Xi\{x\}=1].\ee

\ \

For a set $E\subset \R^d$, let $\Vol_d(E)$ denote the $d$-dimensional volume
of $E$. For $u\in(0,\infty)$, we let $Q_u\subset\R^d$ be the cube centered at the origin having $\Vol_d(Q_u)=u$.

\ \

\noindent{(vi) \bf Non-degeneracy with respect to $\P^{\b \Psi}$.}  
  Say that $\xi$ satisfies non-degeneracy with respect to $\P^{\b \Psi}$ if there exists  $r \in (0, \infty)$
and $b_0:=b_0(r) \in (0, \infty)$
such that  given $\P^{\b \Psi} \cap Q_r^c$, the sum
$ \sum_{x \in \P^{\b \Psi} \cap Q_{t}} \xi(x, \P^{\b \Psi} )$ has expected variability bounded below by
$b_0$,  uniformly in $t \in [r, \infty)$.  In other words, we have
\be \label{assum2}
\inf_{t \in [r, \infty)} \E\Var [ \sum_{x \in \P^{\b \Psi} \cap Q_t } \xi(x, \P^{\b \Psi} )  |
\  \P^{\b \Psi} \cap Q_r^c ] \geq b_0.
\ee

\ \

\noindent As shown in Section~\ref{applications}, functionals of interest often satisfy \eqref{assum2}. There is nothing special about
using cubes $Q_r$ in \eqref{assum2} and, as can be seen from the proofs,
$Q_r$ could be replaced by any compact convex subset of $\R^d$.

\ \
If $f$ and $g$ are two functions satisfying $\liminf_{\la \to \infty} f(\la)/g(\la) > 0$ 
 then
we write $f(\la) = \Omega(g(\la))$.  If, in addition we have $f(\la) = O(g(\la))$ then we write
$f(\la) = \Theta(g(\la))$. 

From the standpoint of applications, it is useful to have a version of \eqref{assum2} for score functions which are not
translation invariant and for input
\be \label{tilP} \tilde{\P}_{\la}^{\b \Psi}:= \P^{\b \Psi} \cap \tilde{S}_\la, \ee  where $\tilde{S}_\la \subset \R^d$ satisfies
$\Vol_d(\tilde{S}_\la) = \Omega(1)$.
Here and elsewhere $\tilde{Q_u} \subset R^d$ denotes a cube with $\Vol_d(\tilde{Q}_u)=u$, but not necessarily centered at the origin. 

\vskip.3cm

\noindent{(vii) \bf  Non-degeneracy  with respect to $\tilde{\P}^{\b \Psi}_\la$.}
Say that $\xi$ satisfies non-degeneracy with respect to $\tilde{\P}^{\b \Psi}_\la$
if there is $r \in (0, \infty)$ and $b_0:=b_0(r) \in (0, \infty)$,
 such that for $\la$ large {there is $\tilde{Q}_r \subset \tilde{S}_\la$ } satisfying
\be \label{noti}
 \E \Var [ \sum_{x \in \tilde{\P}^{\b \Psi}_\la} \xi(x, \tilde{\P}^{\b \Psi}_\la )  |
\  \tilde{\P}^{\b \Psi}_\la \cap \tilde{Q}_r^c ] \geq b_0.
\ee
Given $\rho \in (r, \infty)$, let ${\cal C}(\rho,r, \tilde{S}_\la)$, be a collection of $d$-dimensional volume $r$ cubes
$\tilde{Q}_{i,r}, 1 \leq i \leq n(\rho, r, \tilde{S}_\la),$
which are separated by $\eightrho$ and which satisfy \eqref{noti}.

\vskip.3cm

For all $x$ and $y$ in $\R^d$ we put
\be \label{defcx}
c^\xi(x):= \E \xi(x, \P^{\b \Psi} ) \exp(-
\b \Delta(x, \P^{\b \Psi})),  \ee
and
\be \label{defc2x}
c^\xi(x,y):= c^\xi(x) c^\xi(y) - \E \xi(x, \P^{\b \Psi} \cup \{y\})
\xi(y, \P^{\b \Psi}  \cup \{x \}) \cdot \exp(- \b \Delta(\{x, y \}, \P^{\b \Psi})).
\ee
Put
\be \label{defsig}
\sigma^2(\xi,{\tau}):= c^{\xi^2}(\0) - {\tau}\int_{\R^d} c^\xi(\0, y)dy.
\ee

\subsection{Main results}

The following are our main results.  Applications  follow in Section~\ref{applications}. Our first result
gives conditions under which the Gibbsian functional $W_\la$ has volume order fluctuations.

\begin{theo} \label{main0} 
Assume that
 $\xi$ is translation invariant, exponentially stabilizing in the wide sense \eqref{ExpStab} and
 satisfies the $q$-moment condition \eqref{mom}  for some $q \in (2, \infty)$. Then
 \be \label{Var11}
 \liml \la^{-1} \Var W_\la = {\tau} \sigma^2(\xi, {\tau}) \in [0, \infty).
 \ee
  If, in addition, $\xi$ satisfies non-degeneracy \eqref{assum2}, then $\sigma^2(\xi, {\tau}) > 0$.

\end{theo}



Recall that the Kolmogorov distance between the distributions of random variables $X_1$ and $X_2$ is defined as
$$d_K(X_1,X_2):=\sup_{t\in\RR}|\PP[X_1\le t]-\PP[X_2\le t]|.$$


\begin{theo} \label{main} 
 Assume that
 $\xi$ is exponentially stabilizing in the wide sense
 \eqref{ExpStab} and satisfies the $q$-moment condition \eqref{mom} for some $q \in (2, \infty)$.
 For all $p \in (2, q)$, put $p_3:= p_3(p):= \min\{p,3\}.$
 Then
  \begin{equation}
d_K\left( \frac{W_\la - \E W_\la}{ \sqrt{\Var W_\la} } ,N(0,1)\right) =
 O(    (\ln \la)^{d(p_3 -1)}  \la ( \Var W_\la)^{- p_3/2} ).\label{main1-01}\end{equation}
  Furthermore, if $\xi$ is translation invariant, satisfies non-degeneracy \eqref{assum2} and the $q$-moment condition \eqref{mom}  for some $q \in (3, \infty)$, then
 \begin{equation}
 d_K\left(  \frac{W_\la - \E W_\la}{\sqrt{\Var W_\la} } ,N(0,1)\right) =
 O( (\ln \la)^{2d} \la^{-1/2} ) \label{main1-02}\end{equation}
 and therefore as $\la \to \infty$
 $$
 \la^{-1/2}(W_\la - \E W_\la) \tod N(0, \tau \sigma^2(\xi, \tau)).$$
 \end{theo}
\noindent{\bf Remarks.} (i) (Theorem \ref{main0}.)  The proof of volume order variance asymptotics is indirect. We first
show that $\Var W_\la$ is of volume order up to a logarithmic term (Lemma \ref{varlb}).
Putting $\hat{W}_\la:= \sum_{x \in {\H}_\la^{\b \Psi}} \xi(x,  {\H}^{\b \Psi} \setminusx)$ we then show in Lemma \ref{varld} the dichotomy
that either $\Var \hat{W}_\la = \Omega(\la)$ or $\Var \hat{W}_\la = O(\la^{(d-1)/d}).$   Closeness of $\Var W_\la$ and
$\Var \hat{W}_\la$, as shown in Lemma \ref{var630}, completes the argument, whose full details are in Section 3.
Under condition \eqref{assum2} we obtain volume order variance asymptotics when $\P^{\b \Psi}$ is replaced by a homogeneous Poisson point process, which
is of independent interest.
Verifying condition \eqref{assum2} for Gibbsian input is comparable to verifying
the non-degeneracy conditions of Theorem 2.1 of \cite{PY1} or Theorem 1.2 of \cite{ERS}

\vskip.3cm

\noindent(ii) (Theorem \ref{main}.)  Theorem 2.3 of \cite{SY}
shows the rate of convergence $O( (\ln \la)^{3d} \la^{-1/2})$ in \eqref{main1-02}. However this result
assumes that $\Var W_\la = \Theta(\la)$, which may not always hold, particularly when the scaling is not volume order. 
Theorem \ref{main} contains no such  assumption.
Theorem \ref{main} extends Corollary 3.1 of \cite{BX1} to Gibbsian input.
We do not take up the question of laws of large numbers  for $W_\la$ as this
is  addressed in \cite{SY}.

\vskip.3cm

\noindent(iii) (Point processes with marks.)  Let $({\Ed}, {\cal F}_{{\Ed}}, \mu_{{\Ed}})$  be a probability space (the mark space) and consider
the marked reference Poisson point process 
 $ \{(x,a); x \in  \tilde{\P}_\tau, a \in {\Ed} \}$ in the space $\R^d \times {\Ed}$ with law given by the product measure of the law of  $\tilde{\P}_\tau$  and $\mu_{{\Ed}}$.   Then the proofs of Theorems \ref{main0} and \ref{main} go through in this setting, where it is understood that
 in the algorithmic construction  the process
 ${\P}_{\la}^{\b \Psi}$ inherits the marks
 from $\tilde{\P}_\tau$ and where  the
cubes $Q_r$ in condition \eqref{assum2} are replaced with cylinders $C_r:= Q_r \times {\Ed}$. This generalization is used in Section 2.5 to deduce central limit theorems for spatial birth-growth models with Gibbsian input.



\vskip.3cm
Next we consider the analog of $W_\la$ on input $\tilde{\H}_\la^{\b \Psi}$ defined at \eqref{tilP}, namely
$$
\tilde{W}_\la:= \sum_{x \in \tilde{\H}_\la^{\b \Psi} } \xi(x,  \tilde{\H}_\la^{\b \Psi} \setminus \{x\}).
$$
Say that $\xi$ satisfies the $q$-moment condition with respect to
$\tilde{\H}_\la^{\b \Psi}$ if
\be\label{mom3} \sup_{\la \in [1, \infty) } \sup_{x \in \tilde{S}_\la} \| \xi(x,\tilde{\H}_\la^{\b \Psi} \cup \{x\} )\|_q < \infty. \
\ee
The following result does not assume that $\xi$ is translation invariant.

\begin{theo} \label{main-noti} 
  Assume that
 $\xi$ is exponentially stabilizing in the wide sense
 \eqref{ExpStab} and satisfies the $q$-moment condition \eqref{mom3} for some $q \in (2, \infty)$.
For all $p \in (2, q)$, put $p_3:= p_3(p):= \min\{p,3\}.$  Then
  \begin{equation}
d_K\left(\frac{ \tilde{W}_\la - \E \tilde{W}_\la }{ \sqrt{\Var \tilde{W}_\la} },N(0,1)\right)=
 O\left(    (\ln \la)^{d(p_3 -1)} \Vol(\tilde{S}_\la)   ( \Var \tilde{W}_\la)^{- p_3/2} \right).\label{main1-01-noti}\end{equation}
 Furthermore, if $\xi$ satisfies  non-degeneracy \eqref{noti}
 and $\rho \in (c \ln \la, \infty)$, $c$ large,
  then
  \be \label{VartW}
 \Var \tilde{W}_\la \geq c^{-1}b_0 n(\rho,r, \tilde{S}_\la).\ee  If $q \in (3, \infty)$ we thus have
 \begin{equation}
 d_K\left(\frac{ \tilde{W}_\la - \E \tilde{W}_\la }{ \sqrt{\Var \tilde{W}_\la} },N(0,1)\right) =
 O\left( (\ln \la)^{2d}  \Vol(\tilde{S}_\la) n(\rho,r, \tilde{S}_\la)^{-3/2} \right). \label{main1-02-noti}\end{equation}
 \end{theo}

\noindent{\bf Remark.}  The bound \eqref{VartW} shows volume order growth for $\Var \tilde{W}_\la$, but only up to the
 logarithmic factor  $(\ln \la)^d$. 
When $\xi$ is translation invariant we are able to remove this factor, as described in Remark (i) following Theorem \ref{main}.
However for non-translation invariant $\xi$, we are unable to remove
the logarithmic factor.  Consequently, the  bound  \eqref{main1-02} is smaller than the bound  \eqref{main1-02-noti} by a factor $((\ln\la)^d)^{3/2}$.


\vskip.5cm

\section{ Applications}\label{applications}

\allco


We deduce variance asymptotics and central limit theorems for six well-studied functionals in geometric
probability.  
Save for some special cases as noted below, the limit theory for these functionals has, up to now, been largely confined
to Poisson or binomial input.
Our examples are not exhaustive.  For example, there is
scope for treating the limit theory of coverage processes on Gibbsian input, and, more generally,
the limit theory of functionals of germ-grain models,
with germs given by the realization of ${\cal P}^{\b \Psi}$. One could also treat the limit theory of
functionals arising in percolation and
 nucleation models having Gibbsian input, extending \cite{LP} and
\cite{HQR}, respectively.


\subsection{Clique counts in random geometric graphs}

Let $\X \subset \R^d$ be locally finite and put $s \in (0, \infty)$.  The geometric graph on $\X$, here denoted
$GG_s(\X)$, is obtained by connecting points $x, y \in \X$ with an edge whenever $|x-y| \leq s$.  If there is a subset
$S:= S(s, k)$ of $\X$ of size $k + 1$ with all points of $S$ within a distance $s$ of each other, then the $k$ simplex
formed by $S$ has edges in $GG_s(\X)$.  The Vietoris-Rips complex ${\cal R}^s(\X)$, or Rips complex, is the simplicial complex
arising as the union of of all $k$-simplices $S(s, k) \subset GG_s(\X)$.  The Vietoris-Rips complex and the closely related Cech complex
(which has  a simplex for every finite subset of balls in $GG_s(\X)$ with non-empty intersection) are used to model the topology
of ad hoc sensor and wireless networks and they are also useful in the statistical analysis of high-dimensional data sets.
Note that ${\cal C}^s_k(\X)$  is the number of cliques of order $k + 1$ in $GG_s(\X)$.
 For $\X$ random, the
number ${\cal C}^s_k(\X)$ of $k$-simplices in $GG_s(\X)$ is of theoretical and applied interest (see e.g. \cite{Pe}).    The limit theory for ${\cal C}^s_k(\X)$ is
well understood when $\X$ is Poisson or binomial input on $\R^d$ \cite{Pe} or on a manifold \cite{PY6}.  We are unaware of
limit theory for ${\cal C}^s_k(\cdot)$ on Gibbsian input. For all $k = 1,2,...$ and all $s \in (0, \infty)$ let
$\xi_k(x, \X):= \xi_k^{(s)}(x, \X)$ be $(k + 1)^{-1}$ times the number of $k$-simplices in ${\cal R}^s(\X)$ containing the vertex $x$.

\newpage

\begin{theo} \label{RGG}
For all $k = 1,2,...$ and all $s \in (0, \infty)$
we have
$$
\liml \la^{-1} \Var [{\cal C}^s_k( {\H}_\la^{\b \Psi}) ] = \tau \sigma^2( \xi_k, \tau) > 0,
$$
and
$$
d_K\left(   \frac{   {\cal C}^s_k( {\H}_\la^{\b \Psi}) - \E {\cal C}^s_k( {\H}_\la^{\b \Psi})}  { \sqrt{ \Var [{\cal C}^s_k( {\H}_\la^{\b \Psi}) ]} }, N(0,1) \right) = O(( \ln \la)^{2d} \la^{-1/2} ).
$$
\end{theo}
\noindent {\em Proof}.  We have  ${\cal C}^s_k({\H}_\la^{\b \Psi}) = \sum_{x \in {\H}_\la^{\b \Psi}} {\xi}_k(x, {\H}_\la^{\b \Psi}).$  It suffices to show that $\xi_k$ satisfies the conditions of Theorems  \ref{main0} and \ref{main}.  Given $x \in \R^d$ and $k = 1,2,...$ we note that $\xi_k(x, {\H}_\la^{\b \Psi})$ is generously bounded
by $ \left(  \sum_{X_i \in {\H}_\la^{\b \Psi} } {\bf{1}} (|x - X_i | \leq s) \right)^k $, which is in turn bounded by the $k$th power of a Poisson random variable
with parameter $ \tau \Vol_{d}(B_s(x)) $.  Since all moments of Poisson random variables are finite, it follows that $\xi_k$ satisfies the moment
condition \eqref{mom}  for all $q \in (1, \infty)$.  Clearly $\xi_k$ is translation invariant and exponentially stabilizing
with stabilization radius equal to $s$.  It remains to show that $\xi_k$ satisfies non-degeneracy \eqref{assum2}.  With $s$ fixed, put $r := (3s)^d.$
Let $E_1$ be the event that ${\H}_\la^{\b \Psi}$ puts $k + 1$ points in $Q_{s^d}$ and no points in $Q_r \setminus Q_{s^d}$.  On the event $E_1$ we
have $\sum_{x \in \H_\la^{\b \Psi} \cap Q_r } \xi_k^{(s)}(x, \H_\la^{\b \Psi}) = 1$.  On the other hand, if $E_2$ is the event that ${\H}_\la^{\b \Psi}$ puts fewer than $k + 1$ points in $Q_{s^d}$ and no points in $Q_r \setminus Q_{s^d}$ then $\sum_{x \in \H_\la^{\b \Psi} \cap Q_r } \xi_k^{(s)}(x, \H_\la^{\b \Psi}) = 0$.  Events $E_1$ and $E_2$ have strictly positive probability and give rise to different values of $\sum_{x \in \H_\la^{\b \Psi} \cap Q_r } \xi_k^{(s)}(x, \H_\la^{\b \Psi})$, regardless of the point configurations $\H_\la^{\b \Psi} \cap Q_r^c.$   This shows \eqref{assum2} and concludes the
proof.  \qed


\subsection{Functionals of Euclidean graphs}\label{Euclideangraph}

Many functionals of Euclidean graphs on Gibbsian input satisfy
\eqref{Var11} and \eqref{main1-01}, as shown in \cite{SY}.  However \cite{SY} left open the question of showing variance lower bounds, which
is essential to showing that  \eqref{main1-01} is meaningful.     
We now redress this and assert  that the functionals in \cite{SY}  satisfy  non-degeneracy \eqref{assum2}, and thus
$\sigma^2(\xi, {\tau}) > 0$.  We illustrate this for select functionals in \cite{SY}, leaving it to the reader to verify this assertion
for the remaining functionals, namely those arising in  random sequential adsorption, component counts in random geometric graphs, and Gibbsian loss networks.

\vskip.3cm

\noindent{\bf (i) $k$-nearest neighbors graph.} The $k$-nearest neighbors
(undirected) graph on the vertex set $\X$, denoted $NG(\X)$, is
the graph obtained by including $\{x,y\}$ as an edge
whenever $y$ is one of the $k$ points nearest to $x$ and/or $x$ is
one of the $k$ points nearest to $y$.  The $k$-nearest neighbors
(directed) graph on $\X$, denoted $NG'(\X)$, is obtained by
placing a directed edge between each point and its $k$-nearest
neighbors. In case $\X = \{ x \}$ is a singleton, $x$ has no
nearest neighbor and the {\it nearest neighbor} distance for $x$
is set by convention to $0$.

{\em Total edge length of $k$-nearest neighbors graph.} 
Given $x \in \R^d$ and a locally finite point set $\X \subset \R^d$,
the {\em nearest neighbors length functional} $\xi_{NG}(x,\X)$ is
 one half the sum of the edge lengths of edges in $NG(\X \cup \{x\} )$ which
are incident to $x$.
The total edge length of $NG({\cal P}^{\b \Psi} \cap Q_{\la} )$ is given by
$$W_\la:= \sum_{x \in {\H}_\la^{\b \Psi}} \xi_{NG}(x,  {\H}_\la^{\b \Psi} \setminusx ).$$

Theorem 5.2 in \cite{SY} shows that $W_\la$ satisfies the rate of convergence to the normal
at \eqref{main1-01}. This follows since $\xi_{NG}$  is translation invariant,  exponentially
stabilizing in the wide sense, and  satisfies the moment condition (\ref{mom}) for all $q \in (2, \infty)$.
However that theorem leaves open the question of variance lower bounds
for $\Var W_\la$ and thus the rate of convergence is possibly useless.  The next result
resolves this question {and also gives a slightly better bound than that in \cite{SY}}.

\begin{theo} \label{NNG} 
We have $\liml \la^{-1} \Var W_\la = \tau \sigma^2(\xi_{NG}, \tau) > 0$
and
$$d_K\left(\frac {W_\la-\E W_\la}{\sqrt{\Var W_\la}},N(0,1)\right)=O((\ln \la)^{2d} \la^{-1/2}).$$
 \end{theo}

\noindent{\em Proof.}  We need only  show that non-degeneracy \eqref{assum2} holds
and then apply Theorem  \ref{main0} and \eqref{main1-02}. 
We do this
by modifying the proof of Lemma 6.3 of \cite{PY1}. This goes as follows.
Let $C_0 := Q_{1}$,  the unit cube centered at the origin.
The annulus $Q_{4^d} \setminus C_0$ will be called the moat; notice that $Q_{4^d}$ has
edge length $4$.
Partition the annulus  $Q_{6^d} \setminus Q_{4^d}$ into a
finite collection $\cal U$ of unit cubes. 
Now define the following events. Let
$E_2$ be the event that there are no points in ${\H}_\la^{\b \Psi}$ in the moat
and there are at least $k + 1$ points in each of the unit subcubes in
$\cal U$. Let $E_1$ be the intersection of $E_2$ and the event
that there is $1$ point in $C_0$; let $E_0$ be the intersection of
$E_2$ and the event that there are no points in $C_0$.
Then $E_0$ and $E_1$ have strictly positive probability. Put $Q_r := Q_{6^d}$, i.e., put
$r = 6^d$.

Given any configuration
$\P^{\b \Psi} \cap Q_r^c$, then conditional on the event that $E_0$ occurs,
the sum
$$\sum_{x \in {\H}_\la^{\b \Psi} \cap Q_r} \xi_{NG}(x,  {\H}_\la^{\b \Psi} \setminusx )$$
is strictly less than the same sum, conditional on the event $E_1$. This is because on the event $E_1$
there are $k$ additional edges crossing the moat, each of length at least $3$.

Thus $E_0$ and $E_1$ are events with strictly positive probability
which give rise to values of $\sum_{x \in {\H}_\la^{\b \Psi} \cap Q_r} \xi_{NG}(x,  {\H}_\la^{\b \Psi} \setminus \{x\})$ which differ
by at least $3k$, a fixed amount.  This demonstrates
non-degeneracy \eqref{assum2}.   \qed

\vskip.3cm

\noindent {\bf (ii) Gibbs-Voronoi tessellations}. Given $\X \subset \R^d$ and
$x \in \X$, the set of points in $\R^d$ closer to $x$ than to any
other point of $\X$ is the interior of a possibly unbounded convex
polyhedral cell $C(x,\X)$. The Voronoi tessellation induced by $\X$
is the collection of  cells $C(x,\X), x \in \X$. When $\X$ is the
realization of the Poisson point set $\P_\tau$, this generates
the Poisson-Voronoi tessellation of $\R^d$.
Here, given the Gibbs point process $\P^{\b \Psi}$, consider the Voronoi
tessellation
of this process, sometimes called the Ord process \cite{MW}.  

{\em Total edge length of  Gibbs-Voronoi tessellations.}
Given $\X \subset \R^2$, let $\xi_{ {\rm Vor}}(x, \X)$
denote one half the total edge length of the {\it finite} length
edges in the cell $C(x,\X \cup \{x\})$ (thus we do not take infinite
edges into account).
The total edge length of the Voronoi graph on $\P^{\b \Psi}$ is
given by
$$W_\la:= \sum_{x \in {\H}_\la^{\b \Psi}} \xi_{ \rm{Vor}}(x,  {\H}_\la^{\b \Psi} \setminusx).$$
It may be shown \cite{SY} that $\xi_{\rm{Vor}}$ is exponentially
stabilizing in the wide sense \eqref{ExpStab}, that it satisfies the moment condition (\ref{mom}) 
 for $q \in (2, \infty)$, and, as in Theorem 5.4 of \cite{SY} that $W_\la$ satisfies the rate of convergence to the normal
as in \eqref{main1-01}.

However that theorem leaves open the question of variance lower bounds
for $\Var W_\la$ and thus the rate of convergence is possibly useless.  The next result
resolves this question and gives a better rate than that in \cite{SY}.

\begin{theo} \label{Vor} 
We have $\liml \la^{-1} \Var W_\la = \tau \sigma^2(\xi_{\rm{Vor}}, \tau)  > 0$ and
$$d_K\left(\frac {W_\la-\E W_\la}{\sqrt{\Var W_\la}},N(0,1)\right)=O((\ln \la)^{2d} \la^{-1/2}).$$
\end{theo}

 \noindent{\em Proof.}  We need only  show that non-degeneracy \eqref{assum2} is satisfied
and then apply Theorem  \ref{main0} and \eqref{main1-02}.  We do this
by modifying the proof of Lemma 8.2 of \cite{PY1}. This goes as follows.

 Consider the construction used in the proof of Theorem \ref{NNG}.  Let $E_2$ be the event that there are no points of
${\H}_\la^{\b \Psi}$ in the moat and there is at least one point in each of
the subcubes in $\cal U$. Fix $\eps$ small ($< 1/100$). Choose
points $x_1,x_2,x_3 \in \R^2$ forming an  equilateral triangle of
side-length $1/2$, centered at the origin. 
Let $A_0$ be
the intersection of $E_2$ and the event that there is exactly one
point in each of 
$B_\eps(x_i)$, and
the event that there is no other point in
$C_0{\setminus}\left(\cup_{i=1}^3B_\eps(x_i)\right)$, except for a point $z$ in
the ball $B_{\eps\delta}(0)$,
where $\delta \in (0,1)$ will be chosen shortly.
 Let $A_1$ be the intersection of $E_2$, the event that there is exactly one point in each of $B_{\eps\delta}(\delta x_i)$,
 and the event that there is no other point in
$C_0{\setminus}\left(\cup_{i=1}^3B_{\delta\eps}(\delta x_i)\right)$, except for the point $z$ in
the ball $B_{\eps\delta}(0)$.


On the event $A_0$, the presence of $z$ near the origin
leads to three edges, namely the edges of a
(nearly equilateral) triangular
cell $T$ around the origin. It removes the parts of the three edges of the
Voronoi graph (on all points except $z$)  which  intersect $T$.
The difference between the sum of the lengths of the added edges and
the sum of the
lengths of the three removed edges exceeds some fixed positive
number $\alpha$ (the reason is this: given an
equilateral triangle $T$,
and a point $P$ inside it, the sum of the lengths of the three edges
joining $P$ to the vertices of $T$ is strictly less than the perimeter of
$T$ since the length of each of the
three edges is less than the common length of
the side of $T$. If $T$ is nearly equilateral (our
case) this is still true).

On the other hand, on the event $A_1$, the presence of $z$ cannot increase the total edge length by more than the
total edge length of triangular cell around the origin, and
 this increase is bounded by a constant multiple of
$\delta$, which is less than $\alpha$
if $\delta$ is small enough. Thus if $\delta$ is small enough,
the events $A_0$ and $A_1$ give rise to values of $\sum_{x \in {\H}_\la^{\b \Psi} \cap Q_r} \xi_{\rm{Vor}}(x,  {\H}_\la^{\b \Psi} \setminusx)$
which differ by at least some fixed amount.   This demonstrates
non-degeneracy \eqref{assum2}.
$\qed$


\subsection{Insurance models}




The modeling of insurance claims has been of considerable interest in the literature. The thrust of the modeling is to set up a claim process $\{N_t,\ t\ge 0\}$  to
record the number and time of claims and a sequence of random variables $\{X_i,\ i\ge 1\}$ representing the claim sizes. The aggregate claim size by time $t$ can then be represented
as $S_t=\sum_{i=1}^{N_t}X_i$. Most {of the literature} assumes that $\{X_i,\ i\ge 1\}$ are independent and identically distributed random variables, and are independent of the claim
process $\{N_t,\ t\ge 0\}$ \cite{Embrechts}. When $\{N_t,\ t\ge 0\}$ is a Poisson process, the process $\{S_t, \ t\ge 0\}$ becomes a compound Poisson process and
is also known as
the Cram\'er--Lundberg model (\cite{Embrechts}, p.~22). Significant effort has been devoted to generalize the model so that it represents real situations more closely, e.g.,  making the claim process a more
general counting process such as a renewal process, a negative binomial process, or a stationary point process \cite{Rolski99}. To address the
interdependence of claim sizes, \cite{BX1} introduces a
strictly stationary process $\{Y_t,\ t\ge 0\}$ representing a random environment of the claims and
a simple point process $H$ on $[0,T]\times\N$
recording the times and sizes of clusters of claims. The total claim amount $X_a$
for $a = (t,n)$ is assumed to be the sum of $n$ independent and identically distributed random variables with distribution determined by
the value of $Y_t$. Assuming that $\{Y_t\}$ is independent of $H$ and both $\{Y_t\}$ and $H$ are locally dependent with a `uniform dependence radius $h_0$' such that for all $0<t_1<t_2<\infty$,
$Y\vert_{[t_1,t_2]}$ is independent of $Y\vert_{\R^+{\setminus}(t_1-h_0,t_2+h_0)}$ and
$H\vert_{[t_1,t_2]\times\N}$ is independent of $H\vert_{(\R^+{\setminus}(t_1-h_0,t_2+h_0))\times\N}$, \cite{BX1} proves that the aggregate
claim size $W_T:=\int_{a=(t,n):\ t\le T}X_aH(da)$, when standardized, can be approximated in distribution by the standard normal with
an approximation error of order $O(T^{-1/2})$.

{In disastrous events,  insurance claims may involve dependence amongst the time, size and environment of the claims.
In applications, local dependence with a uniform dependence radius may be violated.}
In this subsection, we aim to address these issues.
 To this end, let the time and spatial location of claims of insurances
be represented by $\P^{\b \Psi}$, a Gibbs point process in $\R^+ \times \R^d$.
In practice, we have $d \in \{2, 3\}$ and the space is typically restricted to a {compact convex}  set $\D\subset \R^d$ with $\Vol_d(\D) > 0$.
Consequently, we set
 $\tilde{\P}_T^{\b \Psi}:=\P^{\b \Psi}\vert_{[0,T]\times\D}$.
Let $\xi((t,\bs), \tilde{\P}_T^{\b \Psi})$  be the value of
the claim at $(t,\bs)$ with $t\in \R^+$ and $\bs\in \R^d$. 
The aggregate claim size
 in the time interval $[0,T]$ is $\tilde{W}_T:=\int_{[0,T]\times \D} \xi((t,\bs), \tilde{\P}_T^{\b \Psi} )\tilde{\P}_T^{\b \Psi}(dt,d\bs).$
 The proof of the next result makes use of Lemma~\ref{newvar} and is thus deferred to Section~\ref{variancemomentbounds}.


\begin{theo} \label{Insurance} 
 Assume that
 $\xi$ is exponentially stabilizing in the wide sense
 \eqref{ExpStab}, translation invariant in the time coordinate $t$, and satisfies the $q$-moment condition \eqref{mom3}
for some $q \in (3, \infty)$.  If there exists an $\epsilon > 0$ such that for all large $T$ there is an
interval $I \subset ( \epsilon T, (1 - \epsilon) T )$ of length $\Theta(1)$, such that
 the conditional distribution $ \tilde{W}_T | \tilde{\P}_T^{\b \Psi} \cap \{([0,T] \setminus I) \times \D \}$ is non-degenerate,
then
 $$ d_K\left(  \frac{\tilde{W}_T - \E \tilde{W}_T}{\sqrt{\Var \tilde{W}_T} } ,N(0,1)\right) =
 O( (\ln T)^{3.5} T^{-1/2} ). $$
 \end{theo}


\begin{coro} 
Assume that the distribution of
 $\xi((t,\bs), \tilde{\P}_T^{\b \Psi})$ is determined by the  $k$-nearest neighbors  of $(t,\bs)$
and satisfies the $q$-moment condition \eqref{mom3}
for some $q \in (3, \infty)$.  If  there exists an $\epsilon > 0$  such that for all large $T$ there is an
interval $I \subset ( \epsilon T, (1 - \epsilon) T )$ of length $\Theta(1)$, such that
 the conditional distribution $ \tilde{W}_T | \tilde{\P}_T^{\b \Psi} \cap \{([0,T] \setminus I) \times \D \}$ is non-degenerate,
    then
 $$ d_K\left(  \frac{\tilde{W}_T - \E \tilde{W}_T}{\sqrt{\Var \tilde{W}_T} } ,N(0,1)\right) =
 O( (\ln T)^{3.5} T^{-1/2} ). $$
\end{coro}

\noindent{\em Proof.} Using the argument of Section~\ref{Euclideangraph} (i), one can easily verify that $\xi$ satisfies all the conditions of Theorem~\ref{Insurance}, hence the conclusion follows.\qed

\subsection{Maximal points of Gibbsian samples}
Let $K:= [0, \infty)^d$. Given $\X \subset \R^d$ locally finite, $x \in
\X$ is called $K$-maximal, or simply maximal if $(K \oplus x) \cap
\X =\{ x\}$.   A point $x= (x_1,...,x_d) \in
\X$ is maximal if there is no other point $(z_1,...,z_d) \in \X$
with $z_i \geq x_i$ for all $1 \leq i \leq d$. The  maximal layer
$m_K(\X)$ is the collection of maximal points in $\X$.  Let
$M_K(\X):= \text{card} (m_K(\X))$.

Consider the region
$$
A:= \{(v,w): \ v \in D, 0 \leq w \leq F(v)\}
$$
where $F: D \to \R$ has continuous partials $F_i, 1 \leq i \leq d
-1,$ bounded away from zero and negative infinity, $D \subset
[0,1]^{d-1},$ and $|F| \leq 1$. We are interested in showing asymptotic normality for
$M_K( [ \la^{-1/d} {\H}_\la^{\b \Psi}  \oplus (1/2,...,1/2) ]  \cap A)$, with ${\H}_\la^{\b \Psi}$ as in \eqref{Gibbs}.
Maximal points are invariant with respect to scaling and translations and it suffices to prove a central limit
theorem for $M_K({\H}^{\b \Psi} \cap \la^{1/d} A).$

The asymptotic behavior and central limit theorem for $M_K(\X)$ with $\X$ either Poisson or binomial input has been
studied in \cite{Dev93,BCHL98,BHLT01,BX,BX1};  the next theorem extends these results to Gibbsian input.

\vskip.5cm

\begin{theo} \label{maxthm} 
We have
 $$
d_K\left( \frac{M_K({\H}^{\b \Psi} \cap \la^{1/d} A ) - \E M_K({\H}^{\b \Psi}  \cap \la^{1/d} A) }{ \sqrt{\Var M_K({\H}^{\b \Psi} \cap \la^{1/d} A  )} },
N(0,1) \right)
=  O\left( (\ln \la)^{(7d-1)/2} \la^{- (d-1)/2d} \right). $$
 \end{theo}

\vskip.5cm

\noindent{\em Proof.} We shall show this is a consequence of Theorem \ref{main-noti} for an appropriate
$\tilde{S}_\la$.
For any subset $E \subset \R^d$ and $\epsilon > 0$ let $E^\epsilon := \{x \in \R^d: \ d(x, E) < \epsilon\},$
where $d(x, E)$ denotes the Euclidean distance between $x$ and the set $E$.
Put $\partial A:= \{(v, F(v)): \ v \in D\}$, 
 $\tilde{S}_\la:= (\la^{1/d} \partial A)^{c \ln \la}$ and in accordance with \eqref{tilP}, we set $\tilde{\H}_\la^{\b \Psi}:= {\H}^{\b \Psi} \cap \tilde{S}_\la$.
Given any $L \in [1, \infty)$ we observe that if $c$ is large then $M_K({\H}^{\b \Psi} \cap \la^{1/d} A ) = M_K(\tilde{\H}_\la^{\b \Psi} \cap \la^{1/d} A   )$
with probability at least $1 - \la^{-L}$.  Since the third moment of $M_K({\H}^{\b \Psi} \cap \la^{1/d} A )$ is bounded by $O(\la^3)$, this is enough to guarantee that $\Var M_K({\H}^{\b \Psi} \cap \la^{1/d} A  )$ and
$\Var M_K(\tilde{\H}_\la^{\b \Psi}\cap \la^{1/d} A )$ have the same asymptotic behavior and thus it is enough to prove Theorem \ref{maxthm}
with ${\H}^{\b \Psi} \cap \la^{1/d} A$ replaced by $\tilde{\H}_\la^{\b \Psi} \cap \la^{1/d} A$.
Put  \bean \zeta(x, \X) := \zeta(x, \X; \la^{1/d}A):= \left\{
\begin{array}{ll} 1  & \text{ if }   ((K \oplus x) \cap \la^{1/d} A) \cap (\X \cup \{x\}) = \{x\},
\\
0  & \text{ otherwise}.
\end{array}
\right. \eean
Notice that $\zeta$ is not translation invariant and that
$$
M_K(\tilde{\H}_\la^{\b \Psi} \cap \la^{1/d} A) = \sum_{x \in \tilde{\H}_\la^{\b \Psi} \cap \la^{1/d} A }  \zeta(x,
\tilde{\H}_\la^{\b \Psi}).
$$
To prove Theorem \ref{maxthm}, it suffices to show that $\zeta$ satisfies exponential stabilization
in the wide sense  \eqref{ExpStab} and apply  Theorem \ref{main-noti}. 

To show exponential stabilization, we argue as follows.
Given $x \in \tilde{S}_\la\cap \la^{1/d} A$, let $D_1(x):= D_1(x, \tilde{\H}_\la^{\b \Psi})$ be the distance between
$x$ and the nearest point in $(K \oplus x) \cap \la^{1/d} A \cap \tilde{\H}_\la^{\b
\Psi}$, if there is such a point; otherwise we let $D_1(x)$ be the {\em maximal} distance  between $x$ and $(K \oplus
x) \cap
\la^{1/d} \partial A$, denoted here by $D(x)$. By the smoothness assumptions on
$\partial A$, it follows that $(K \oplus x) \cap \la^{1/d}A  \cap B_t(x)$ has volume at least $c_1 t^d$ for all
$t \in [0, D(x)]$.
It follows that uniformly in $x \in \tilde{S}_\la\cap \la^{1/d} A$ and $\la \in [1, \infty)$
\be \label{expz} \PP[D_1(x) > t]
\leq \exp(-c_1 t^d), \ 0 \leq t \leq D(x).
 \ee For $t \in (D(x), \infty)$,
this inequality holds trivially and so \eqref{expz} holds for all $t \in (0, \infty)$.


Let $R(x):=R(x, \tilde{\H}_\la^{\b \Psi}) := D_1(x)$.  We claim
that $R:=R(x)$ is a radius of stabilization for $\zeta$ at $x$. Indeed,
if $D_1(x) \in (0, D(x))$, then $x$ is not maximal, and so
$$
\zeta(x, \tilde{\H}_\la^{\b \Psi} \cap B_R(x))= 0$$ and inserting  points $\cal Y$
outside $B_R(x)$ does not modify the score $\zeta$. If $D_1(x) \in
[D(x), \infty)$ then
$$
\zeta(x, \tilde{\H}_\la^{\b \Psi} \cap B_R(x))= 1.$$
{Keeping the realization $\tilde{\H}_\la^{\b \Psi} \cap B_R(x)$ fixed, we notice that inserting points
$\cal Y$ outside $B_R(x)$ does not modify the score $\zeta$, since maximality
of $x$ is preserved.} 
Thus $R(x)$ is a
radius of stabilization for $\zeta$ at $x$ and it decays
exponentially fast, as demonstrated above.


Clearly the moment condition \eqref{mom3} is satisfied since $\zeta$ is bounded by one.
We now show that $\zeta$ satisfies  non-degeneracy
 \eqref{noti} for a large number of cubes of volume at least ${c_2}r$.  We do this for $d = 2$, but
the proof extends to higher dimensions.

Fix $r \in [1, \infty)$.
Let $\tilde{Q}_r \subset \tilde{S}_\la$ be such that $\tilde{Q}_r  \cap \la^{1/d} \partial A \neq \emptyset.$
We also assume that $\la^{1/d} A$ contains only the lower left corner of $\tilde{Q}_r$, but that
$\Vol( \tilde{Q}_r  \cap \la^{1/d} A) \geq c_3 r$.


\begin{wrapfigure}{r}{0.55\textwidth}
  \vspace{-20pt}
  \begin{center} 
   \includegraphics[trim = 55mm 20mm 45mm 55mm, clip,width=0.5\textwidth]{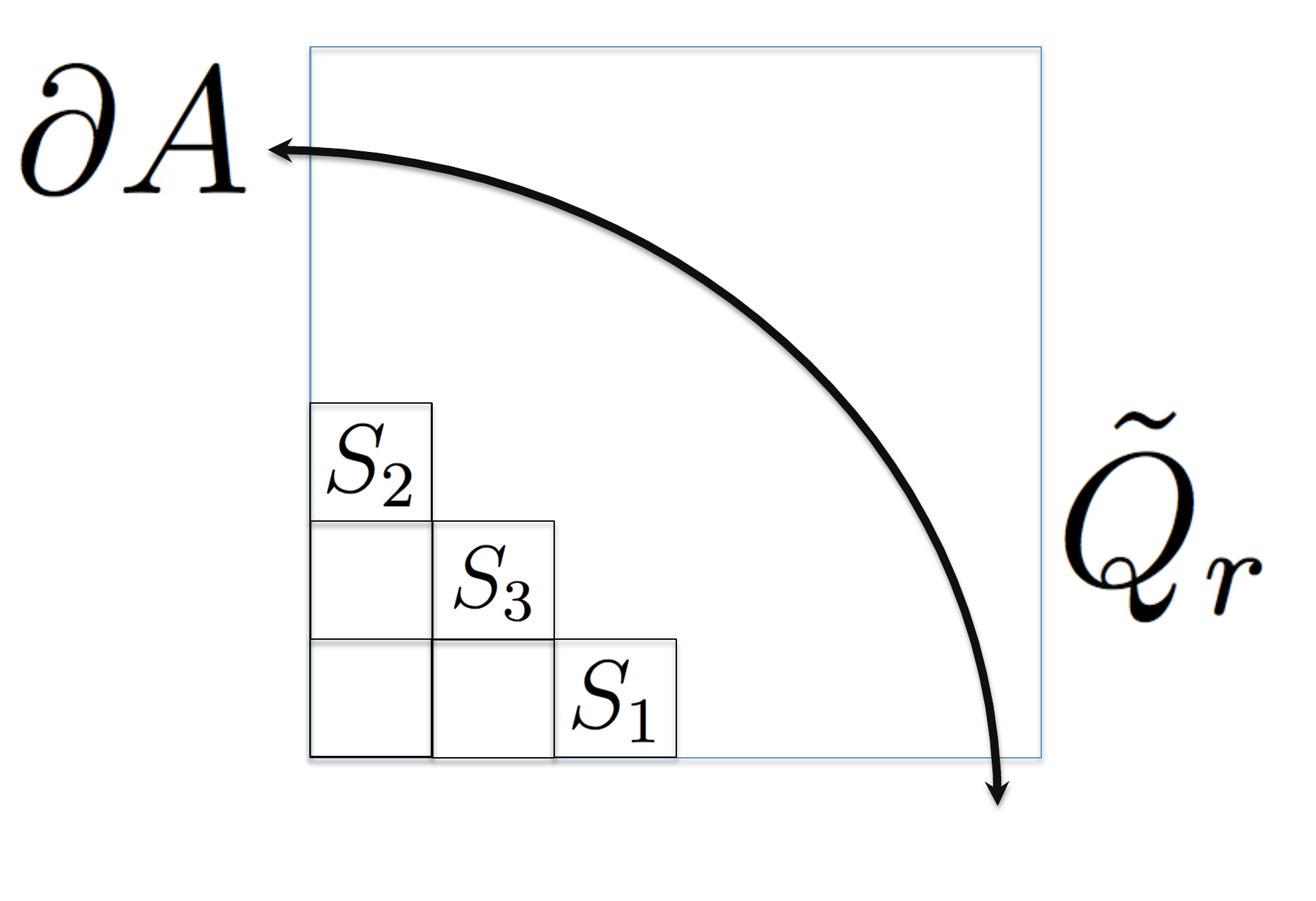}
 \vspace{-20pt}
  \caption{The square $\tilde{Q}_r$ and the subsquares $S_1,S_2, S_3$} %
  \label{figureone}
  \end{center}
  \vspace{-15pt}
 \end{wrapfigure}


Referring to {Figure 1},
we consider the event $E$ that
$\rm{card}( \tilde\H_\la^{\b \Psi} \cap S_1) =
\rm{card}( \tilde\H_\la^{\b \Psi} \cap S_2) = 1$, where $S_1$ and $S_2$ are the squares in {Figure 1}.  Let $E_1$ be the event that $\tilde\H_\la^{\b \Psi}$
puts no  points in $\tilde{Q}_r {\setminus} (S_1\cup S_2)$. Note that $\PP[E\cap E_1]$ is bounded away from zero, uniformly in $\la$. 
On $E\cap E_1$ we have that
$$
\sum_ {x \in \tilde{\H}_\la^{\b \Psi} \cap \tilde{Q}_r   } \zeta(x,  \tilde{\H}_\la^{\b \Psi} )
$$
contributes a value of $2$ to the total sum $\sum_{x \in \tilde{\H}_\la^{\b \Psi}  } \zeta(x,  \tilde{\H}_\la^{\b \Psi} ).
$
Let $E_2$ be the event that $\tilde{\H}_\la^{\b \Psi}$
puts no  points in $\tilde{Q}_r {\setminus} (S_1\cup S_2)$, except for a singleton in the square $S_3$.   Then
$\PP[E\cap E_2]$ is bounded away from zero, uniformly in $\la$.
On $E\cap E_2$ we have that
$$
\sum_{x \in \tilde{\H}_\la^{\b \Psi} \cap \tilde{Q}_r  } \zeta(x,  \tilde{\H}_\la^{\b \Psi} )
$$
contributes a value of $3$ to the total sum $\sum_{x \in \tilde{\H}_\la^{\b \Psi}   } \zeta(x,  \tilde{\H}_\la^{\b \Psi} ).
$
This is true regardless of the configuration ${\tilde{\H}_\la^{\b \Psi}}  \cap \tilde{Q}_r^c$  and so
condition  \eqref{noti} holds.
Since the surface area of $\la^{1/d} \partial A$ is $\Theta( \la^{(d-1)/d} )$,
the number of cubes $\tilde{Q}_r$ having these properties is of order $\Theta((\la^{1/d} / \ln \la)^{d-1})$, whenever $\rho = \Theta(\ln \la)$.
Thus we have  $n(\rho, {r}, \tilde{S}_\la) =
\Theta((\la^{1/d} / \ln \la)^{d-1})$. 

Applying  Theorem \ref{main-noti} we obtain Theorem \ref{maxthm}.
Noting that $\Vol_{d}(\tilde{S}_\la) = \Theta(\la^{(d-1)/d} \ln \la)$, the bound \eqref{main1-02-noti} yields the rate of convergence to the normal
$$
=O\left( ( \ln \la)^{2d} \la^{(d-1)/d} \ln \la \  ( \la^{(d-1)/d} / (\ln \la)^{d-1} )^{-3/2}\right)
= O\left(( \ln \la)^{7d/2 - 1/2} \la^{-(d-1)/2d}\right),
$$
which was to be shown. \qed

\subsection{Spatial birth-growth models}   Consider the following spatial birth-growth model on $\R^d$.
Seeds appear at random locations $X_i \in \R^d$ at {i.i.d.} times
$T_i, \ i = 1,2,...$ according to a
spatial-temporal point process $\P := \{(X_i, T_i) \in \R^d
\times [0, \infty)\}$. 
When a
seed is born, it has initial radius zero and then forms a cell
within $\R^d$ by growing radially in all directions with a constant
speed $v > 0$. Whenever one growing cell touches another, it stops
growing in that direction.
If a seed appears at $X_i$ and if $X_i$
belongs to any of the cells existing at the time $T_i$,
then the seed is discarded.  We assume that the law of $X_i, i \geq 1,$ is independent of the law of
$T_i, i \geq 1.$

Such growth models have received
considerable attention  with mathematical
contributions given in  \cite{CL, CQ,CQa, HQR, Pe1}.
First and second order characteristics for Johnson-Mehl growth models
on homogeneous Poisson points on $\R^d$ are given in \cite{Mo1, Mo2}.
Using the general Theorem  \ref{main}, we may extend many of these results to
growth models with Gibbsian input.  We illustrate with the following theorem, in which $\P$ denotes
a marked Gibbs point process with intensity
measure $m^{\b \psi} \times \mu$, where $m^{\b \psi}$ is the intensity measure of
${\H}^{\b \Psi}$ and $\mu$ is an arbitrary  probability measure on $[0, \infty).$

Given a compact subset $K'$ of $\R^d$, let $N(\P; K')$ be the number of seeds accepted in $K'$.
We shall deduce the following result from Remark (iii) following Theorem \ref{main}. We let $\hat{\H}_\la^{\b \Psi}$ denote
the process of marked points $\{ (X_i, T_i): \  X_i \in {\H}_\la^{\b \Psi}, \ T_i \in [0, \infty) \}$.
Given a marked point set $\X \subset \R^d \times [0, \infty)$, define the score \bean \nu(x, \X) := \left \{
\begin{array}{ll} 1  & \text{ if the seed at $x$ is accepted,}
\\
0  & \text{ otherwise}.
\end{array}
\right. \eean


\begin{theo} \label{sp-timethm} We have
$\lim_{\la \to \infty}  \la^{-1} \Var  N( \hat{\H}_\la^{\b \Psi}; Q_\la ) = \tau \sigma^2( \nu, \tau) > 0$
and
 $$
d_K\left( \frac{ N(\hat{\H}_\la^{\b \Psi}; Q_\la ) - \E N(\hat{\H}_\la^{\b \Psi}; Q_\la ) }
 { \sqrt{\Var N( \hat{\H}_\la^{\b \Psi}; Q_\la )} },N(0,1) \right)
 =
  O\left(  (\ln \la)^{2d} \la^{-1/2} \right).$$
 \end{theo}

\noindent {\em Proof.}  Notice by the definition of $\nu$ we have
 $$N(\hat{\H}_\la^{\b \Psi}; Q_\la ) = \sum_{x \in \hat{\H}_\la^{\b \Psi} \cap Q_\la} \nu(x, \hat{\H}_\la^{\b \Psi}).$$

Let $K$ denote the downward right circular cone with apex at the origin of $\R^d$.  Then
\bean \nu(x, \X) = \left\{
\begin{array}{ll} 1  & \text{ if} \ (K \oplus x) \cap (\X \cup \{x\} ) = x,
\\
0  & \text{ otherwise}.
\end{array}
\right. \eean

We now aim to show that $\nu$ satisfies all the conditions of Theorem \ref{main}.
Clearly $\nu$ is translation invariant in $\R^d$.  The moment condition \eqref{mom} is satisfied,
since $|\nu| \leq 1$.  We claim that $\nu$ satisfies exponential stabilization in the wide sense.  This however follows from the above proof
that $\zeta$ is exponentially stabilizing in the wide sense (the proof is easier now because the boundary of $A$ corresponds to the hyperplane $\R^d$).

We claim that non-degeneracy \eqref{assum2} holds. But this too follows from simple modifications of the proof of non-degeneracy of
$\zeta$.  In fact things are easier, because we need only show  that \eqref{assum2} holds for one cube $Q_r$.  
To this end, the cube $Q_r$ is
now replaced by a space-time cylinder $C_r:= [-r^{1/d}, r^{1/d}]^d \times [0, \infty)$. For simplicity of exposition only, we show non-degeneracy for ${d = 1}$, but the
approach extends to all dimensions.

\begin{wrapfigure}{r}{0.55\textwidth}
  \vspace{-30pt}
  \begin{center} 
   \includegraphics[trim = 65mm 20mm 45mm 15mm, clip,width=0.45\textwidth]{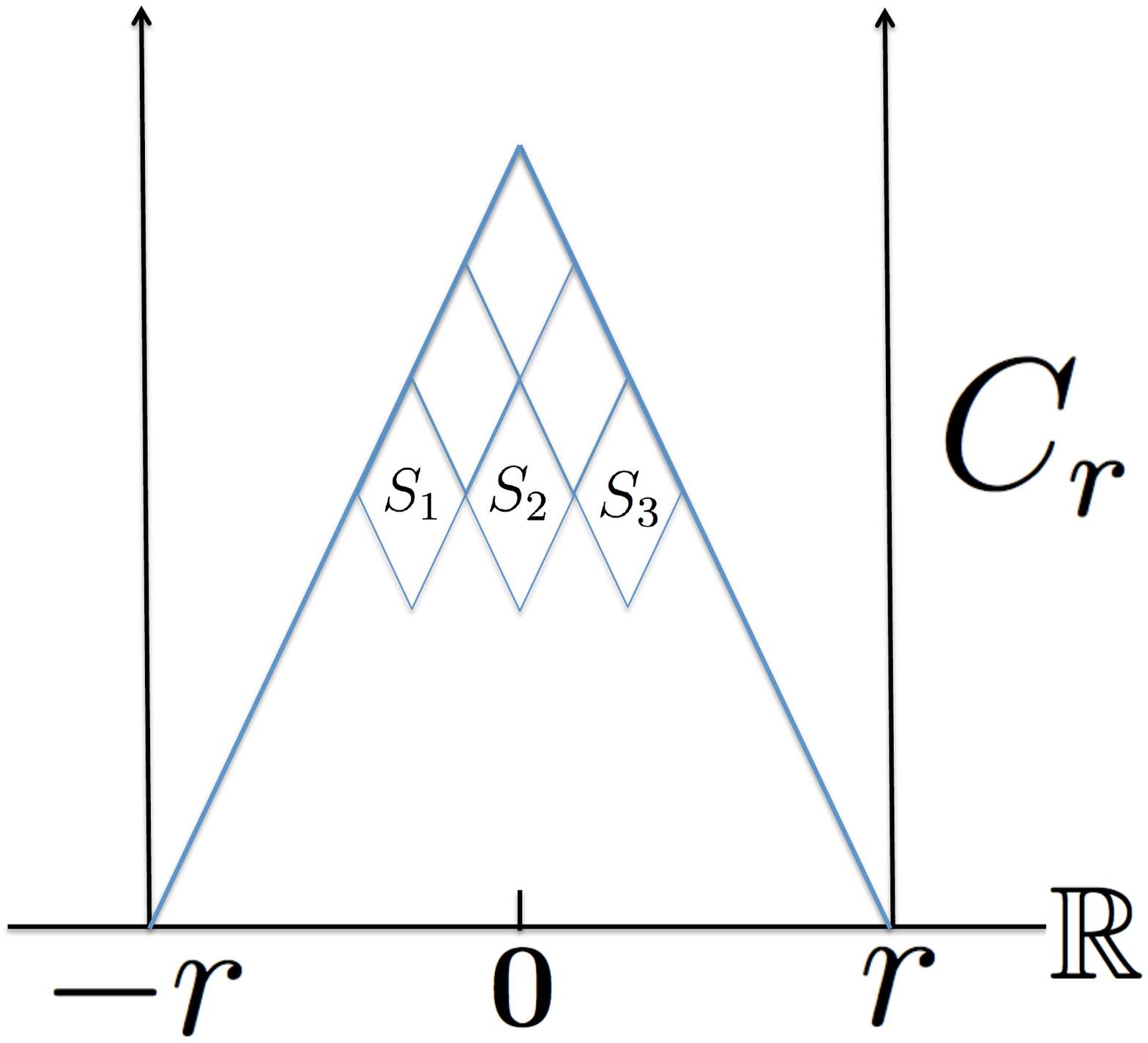}
  \end{center}
    \vspace{-30pt}
  \caption{Space-time cylinder $C_r$}
  \label{figuretwo}
  \vspace{12pt}
\end{wrapfigure}


Referring to {Figure 2}, we consider the event $E$ that $\rm{card}( {\hat{\H}}_\la^{\b \Psi} \cap S_1) =
\rm{card}( {\hat{\H}}_\la^{\b \Psi} \cap S_2) = 1$.  Let $E_1$ be the event that ${\hat{\H}}_\la^{\b \Psi}$
puts no other points in $({[-r, r]} \times [0, 1]){\setminus}(S_1\cup S_2)$ (we don't care about the point configuration in the set
${[-r, r]} \times [1, \infty)$.  Note that $\PP[E\cap E_1]$ is bounded away from zero, uniformly in $\la$.
On $E\cap E_1$ we have that
$$
\sum_{x \in \hat{\H}_\la^{\b \Psi} \cap C_r} \nu(x,  \hat{\H}_\la^{\b \Psi} )
$$
contributes a value of $2$ to the total sum $\sum_{x \in \hat{\H}_\la^{\b \Psi} \cap C_t   } \nu(x,  \hat{\H}_\la^{\b \Psi} ).
$
Let $E_2$ be the event that ${\hat{\H}}_\la^{\b \Psi}$
puts no other points in $({[-r, r]} \times [0, 1]){\setminus}(S_1\cup S_2)$, except for a singleton in the diamond $S_3$.   Then
$\PP[E\cap E_2]$ is bounded away from zero, uniformly in $\la$. On $E\cap E_2$ we have that
$$
\sum_{x \in \hat{\H}_\la^{\b \Psi} \cap C_r } \nu(x,  \hat{\H}_\la^{\b \Psi} )
$$
contributes a value of $3$ to the total sum $\sum_{x \in \hat{\H}_\la^{\b \Psi} \cap C_t  } \nu(x,  \hat{\H}_\la^{\b \Psi} ).
$
This is true regardless of the configuration ${\hat{\H}}_\la^{\b \Psi}  \cap C_r^c$  and so
condition  \eqref{assum2} holds. {Thus}  $\nu$ satisfies all conditions of Theorem~\ref{main} and so
Theorem~\ref{sp-timethm} follows.   \qed


\section{ Auxiliary results }

\allco

Before proving our main theorems we require a few additional results.

 \vskip.5cm

\noindent{(i) \bf Control of spatial dependencies of Gibbs point processes}.
Recall that $\P^{\b \Psi}$ is an admissible point process, i.e.,  $\Psi \in \bf{\Psi}^*$ and $(\tau,\beta) \in {\cal R}^\Psi$.
As shown in the perfect simulation
 techniques of \cite{SY}, the process has spatial dependencies which can be controlled by the size of the so-called
 {\em ancestor clans}.  The ancestor clans are backwards in time oriented percolation clusters,
 where two nodes in space time are linked with a directed edge if one is the ancestor of the other.  The acceptance status of a point
 at $x$ depends on points in the ancestor clan.   As seen at (3.6) of \cite{SY}, the ancestor clans
 have exponentially decaying spatial diameter. Thus, if $A_B^{\b \Psi}(t)$ is the ancestor clan
 {in $\P^{\b \Psi}$ of the set} $B\subset\R^d$ at time $t$, 
  then for all
  $(\tau,\beta) \in {\cal R}^\Psi$, there is a constant $c:=c(\tau,\beta) \in (0, \infty)$
 such that for all $t \in (0, \infty)$, $M \in (0, \infty)$, and $B \subset \R^d$ we have
 \be \label{controlsdep1}
 \PP[ \diam(A_B^{\b \Psi}(t)) \geq M + \text{diam}(B)] \leq c (1+\text{vol}(B))\exp(-M/c).
 \ee
Let $A_{B, \la}^{\b \Psi}$ be the ancestor clan { in $ {\H}_\la^{\b \Psi}$ of the set $B$}. 
Since $\diam(A_{B, \la}^{\b \Psi}(t))\le\diam(A_B^{\b \Psi}(t))$,  the bound \Ref{controlsdep1} also holds for $A_{B, \la}^{\b \Psi}$, i.e., for all $\la \in [1, \infty)$, $B \subset Q_\la$ we have
$$
 \PP[ \diam(A_{B, \la}^{\b \Psi}(t)) \geq M + \text{diam}(B)] \leq c (1+\text{vol}(B))\exp(-M/c).
$$

 Put for all $\rho \in (0, \infty)$ 
$$
d(\rho):= \limsup_{\la \to \infty} \sup_{B \subset Q_\la, \ \text{diam}(B) \leq \rho/2 }  \PP[ \diam(A_{B, \la}^{\b \Psi}) \geq \rho].
$$
Then we have \be \label{diambd} d(\rho) \leq c (1 + (\rho/2)^d v_d) \exp(-\rho/2 c).\ee

\vskip.5cm

\noindent{(ii)  \bf Score functions with deterministic range of dependency.} {
Given the radius of stabilization  $R^\xi(x, {\H}_\la^{\b \Psi})$, let
$D(x, {\H}_\la^{\b \Psi})$ be the diameter of the ancestor clan of the
stabilization ball $B_{R^\xi(x, {\H}_\la^{\b \Psi})}(x)$.
For all $\rho \in (0, \infty)$, consider score functions on points having ancestor clan {diameter} at most $\rho$: 
$$\xi(x,  {\H}_\la^{\b \Psi} \setminusx ; \rho):= \xi(x,  {\H}_\la^{\b \Psi} \setminusx) 
{\bf 1}(D(x, {\H}_\la^{\b \Psi}) \leq  \tworho).
$$
We study the following functional, the analog of $W(\rho)$ on page 704 of \cite{BX1}:
\be \label{Wla}
W_\la(\rho):= \sum_{x \in {\H}_\la^{\b \Psi}} \xi(x,  {\H}_\la^{\b \Psi} \setminusx;
\rho).\ee
When sets $A$ and $B$ are separated by a Euclidean distance greater than $2 \rho$, then
the random variables
$\sum_{x \in {\H}_\la^{\b \Psi} \cap A} \xi(x,  {\H}_\la^{\b \Psi} \setminusx;
\rho)$ and  $\sum_{x \in {\H}_\la^{\b \Psi} \cap B} \xi(x,  {\H}_\la^{\b \Psi} \setminusx;
\rho)$
{depend on disjoint and hence independent portions of the birth and death process $(\varrho(t))_{t \in \R}$ in the construction of
$\P_\la^{\b \Psi}$.}  We make heavy use of this in the proofs of Theorems \ref{main} and \ref{main-noti}.

It is also useful to consider sums of scores with respect to the global point process
${\H}^{\b \Psi}$, namely
$$
\hat{W}_\la:= \sum_{x \in {\H}_\la^{\b \Psi}} \xi(x,  {\H}^{\b \Psi} \setminusx); \ \
\hat{W}_\la(\rho):= \sum_{x \in {\H}_\la^{\b \Psi}} \xi(x,  {\H}^{\b \Psi} \setminusx;
\rho).$$  }

\vskip.5cm

\noindent {(iii) \bf Wide sense stabilization of $\xi$ on ${\H}_\la^{\b \Psi}$}.
If $\xi$ is a stabilizing functional in the wide sense, then  
$$
Q(\rho):=\limsup_{\la \to \infty}\sup_{x\in Q_\la}\PP [
R^{\xi}(x, {\H}_\la^{\b \Psi}) > \rho|{\H}_\la^{\b \Psi}\{x\}=1]
\to 0,
$$
as $\rho \to \infty$. 
If $\xi$ is {\em exponentially stabilizing} in the wide sense \eqref{ExpStab}, then by  \eqref{ExpStabA} there is a constant
${c} \in (0, \infty)$ such that
\be \label{expo}
Q(\rho) \leq {c} \exp(-\rho/{c}).\ee
Notice that for any $\rho \in (0, \infty)$ we have
\begin{eqnarray*}
&&\PP[ D(x, {\H}_\la^{\b \Psi}) \geq \rho|{\H}_\la^{\b \Psi}\{x\}=1 ] \\
&\leq& \PP[ D(x, {\H}_\la^{\b \Psi}) \geq \rho, R^\xi(x, {\H}_\la^{\b \Psi} ) \leq  \rho/2|{\H}_\la^{\b \Psi}\{x\}=1 ]\\
&&+ \PP[  R^\xi(x, {\H}_\la^{\b \Psi} ) \geq  \rho/2 |{\H}_\la^{\b \Psi}\{x\}=1].
\end{eqnarray*}
Bounding the first term on the right hand side by \eqref{diambd} and the second by \eqref{expo},
we obtain  whenever $\rho \in [c' \ln \la, \infty)$ and $c'$ is large that
there is $c_1$ such  that  {$\PP[ D(x, {\H}_\la^{\b \Psi}) \geq \rho|{\H}_\la^{\b \Psi}\{x\}=1 ] \leq {c_1} \exp(-\rho/{c_1})$}  whenever $\rho \in [c' \ln \la, \infty)$.
Thus,   for any $L \in [1, \infty)$,  there is $c$ large enough so that if $\rho \in [c \ln \la, \infty)$,  then
 \be \label{bd2a}
\PP[ \hat{W}_\la \neq \hat{W}_\la(\rho) ] \leq \la^{-L}
\ee 
 and
\be \label{bd2b}
\PP[ W_\la \neq W_\la(\rho) ] \leq \la^{-L}.
\ee

\section{Variance and moment bounds}\label{variancemomentbounds}

\allco

Let $r$ satisfy non-degeneracy \eqref{assum2} and let $\rho \in [r, \infty)$. Find a maximal collection of disjoint cubes
$Q_{i,r}:= Q_{i,r,\rho} \subset Q_\la, i \in I,$ {with $\Vol_{d} Q_{i,r} = r$,}  and which are separated by a distance at least $\eightrho$ and which are at least a distance $2\rho$ from
$\partial Q_\la$. Notice that
$n(\rho, Q_\la):= {\rm{card}}(I) = \lfloor c' \la/\rho^d\rfloor$, $c'$ a constant. 
Let $\F_i$ be the sigma algebra generated by ${\H}^{\b \Psi} \cap Q_{i,r}^c$. More precisely, letting ${\cal B}$ be the class of all locally finite subsets of $\R^d$, define the sigma algebra ${\mathscr B}$ in ${\cal B}$ as the smallest sigma algebra making the mappings $\eta \in{\cal B}\mapsto\card(\eta\cap \Theta),$ for all Borel sets $\Theta\subset\R^d$, measurable (see \cite{Kallenberg83}, page~12). The sigma algebra $\F_i$ is induced by the {mapping
${\H}^{\b \Psi} \mapsto {\H}^{\b \Psi}\cap Q_{i,r}^c$} from ${\cal B}$ to $({\cal B},{\mathscr B})$.


\begin{lemm} \label{momentlemma} Let $q \in [1, \infty)$. If $\xi$ satisfies the moment condition \eqref{mom}
for some $q' \in (q, \infty)$ then  there {are constants $\la_0\in(0,\infty)$ and $c\in(0,\infty)$ such that for all $\la\ge \la_0$ and $\rho \in [1, \infty)$}
\be \label{bdd1}
\max \{ \| W_\la \|_q, \|W_\la(\rho)\|_q \} \leq c \la
\ee
and
$$\sup_{i \in I} \max \{ \| \E[W_\la | {\cal F}_i] \|_q, \|  \E [W_\la(\rho) | {\cal F}_i] \|_q \} \leq c \la.
$$
Identical bounds hold if  $W_\la$ is replaced by $\hat{W}_\la$.
\end{lemm}

\noindent{\em Proof.} Fix $q \in [1, \infty)$. We shall only prove $\| W_\la \|_q \leq c \la$ as the other inequalities follow
similarly.  Put $N:= \text{card}(\P_\la^{\b \Psi})$. Minkowski's inequality gives
$$\|W_\la\|_q \leq \sum_{j = 0}^{\infty} \|\sum_{x \in \P_\la^{\b \Psi}} \xi(x, \P_\la^{\b \Psi} \setminusx  )
{\bf {1}}(\la \tau 2^j \leq N \leq \la \tau 2^{j + 1})\|_q$$
$$ \leq   \sum_{j = 0}^{\infty} \|\sum_{x \in \P_\la^{\b \Psi}, \  N \leq  \la \tau 2^{j + 1} } \xi(x, \P_\la^{\b \Psi}  \setminusx)
{\bf {1}}(N \geq \la \tau 2^{j })\|_q.
$$
Let $s \in (1, \infty)$ be such that $qs < q'$.  Let $1/s + 1/t = 1$, i.e., $s$ and $t$ are conjugate exponents.
H\"older's inequality gives
$$\|W_\la\|_q \leq \sum_{j = 0}^{\infty} \left[ \E ( \sum_{x \in \P_\la^{\b \Psi}, \  N \leq  \la \tau 2^{j + 1} } \xi(x, \P_\la^{\b \Psi} \setminusx  ))^{qs}
\right]^{1/qs} (\PP[N \geq  \la \tau 2^{j}])^{1/qt}.$$
Since $\P_\la^{\b \Psi}$ is Poisson-like, we have that $N$ is stochastically dominated by a Poisson random variable ${\rm Po}(\lambda \tau)$
with parameter $\la \tau$.   Recalling the definition of $w_q$ at \eqref{mom}, we obtain
$$\|W_\la\|_q \leq \sum_{j = 0}^{\infty} \| \sum_{x \in \P_\la^{\b \Psi},  N \leq  \la \tau 2^{j + 1} } \xi(x, \P_\la^{\b \Psi}  \setminusx ) \|_{qs}
 (\PP[{\rm Po}(\lambda \tau)  \geq  \la \tau 2^{j}])^{1/qt}$$
 $$
 \leq 6 \la \tau w_{qs} + \sum_{j = 2}^{\infty} \la \tau 2^{j + 1} w_{qs}
 (\PP[{\rm Po}(\lambda \tau) - \lambda \tau  \geq  \la \tau (2^{j} -1) ])^{1/qt},$$
using Minkowski's inequality another time. For $j \geq 2$, we have that $\PP[{\rm Po}(\lambda \tau) - \lambda \tau  \geq  \la \tau (2^{j} -1)]$  decays
exponentially fast in $2^j$ by standard tail probabilities for the Poisson random variable.
This shows that the infinite sum is $O(\la \tau)$, concluding the proof.   \qed

\  \

We put
$$
\tilde{W}_\la(\rho):= \sum_{x \in \tilde{\H}_\la^{\b \Psi} }\xi(x, \tilde{\H}_\la^{\b \Psi} \setminus \{x\}; \rho).
$$


\begin{lemm} \label{newvar} {Given a set $G \subset \R^d$ we let $\G_G$ (respectively $\tilde{\G}_G$)
 be the sigma algebra generated by ${\H}^{\b \Psi} \cap G$ (respectively $\tilde{{\H}}_\la^{\b \Psi} \cap G$).
Assume that $\xi$ satisfies condition \eqref{ExpStab}.}

\noindent(a)  If $\xi$ satisfies the moment condition \eqref{mom}  for some $q \in (2, \infty)$, then
there exist constants $\la_0$ and $c$ such that for all $\la \in [\la_0, \infty)$, $\rho \in [c \ln\la, \infty)$ and all Borel sets $G\subset\R^d$,
\begin{equation} |\E \Var [ \hat{W}_\la(\rho) |  \G_G]-\E \Var [ \hat{W}_\la |  \G_G] |\le \la^{-1} \label{Xia14082901}\end{equation}
 and
 \begin{equation} |\E \Var [ W_\la(\rho) |  \G_G]-\E \Var [ W_\la |  \G_G] |\le \la^{-1}.\label{Xia14082902}\end{equation}
\noindent
(b) If  $\xi$ satisfies the moment condition \eqref{mom3}  for some $q \in (2, \infty)$ then
{there exist constants $\la_0\in (0, \infty)$ and $c\in (0, \infty)$ such that}
for all $\la \in [\la_0, \infty)$, $\rho \in [c \ln\la, \infty)$ and all Borel sets $G\subset \tilde{S}_\la$,
 $$ |\E \Var [ \tilde{W}_\la(\rho) |  \tilde{\G}_G] -\E \Var [ \tilde{W}_\la |  \tilde{\G}_G]| \le \la^{-1}.$$
 \end{lemm}

\noindent{\em Proof.}  (a) 
Using the generic formula $\Var[X | \A] = \E [X^2 | \A] - (\E[X| \A])^2$, valid for any random variable $X$ and sigma algebra $\A$, we have
$$
\E \Var [ \hat{W}_\la(\rho) |  \G_G] = \E \left[ \E  [\hat{W}_\la^2(\rho) |  \G_G] -  (\E [\hat{W}_\la(\rho) |  \G_G)^2 \right]
$$
and
$$
\E \Var [ \hat{W}_\la |  \G_G] = \E \left[ \E  [\hat{W}_\la^2 |  \G_G] -  (\E [\hat{W}_\la |  \G_G)^2 \right].
$$
If both differences
\be \label{diff1}
| \E [ \E  [\hat{W}_\la^2(\rho) |  \G_G] - \E  [\hat{W}_\la^2 |  \G_G]  ] |
\ee
and
\be \label{diff2}
 | \E [ \E  [\hat{W}_\la(\rho) |  \G_G]^2 - \E  [\hat{W}_\la |  \G_G]^2  ] ] |
\ee
are less than $\la^{-1}/2$ then $\E \Var [\hat{W}_\la(\rho) |  \G_G]$ differs
from $\E \Var [ \hat{W}_\la| \G_G]$ by less than $\la^{-1}$.

Notice that \eqref{diff1}  may be bounded by $(2\la)^{-1}$ since it equals $\E [\hat{W}_\la^2(\rho) -\hat{W}_\la^2]$, which by
H\"older's inequality is bounded by the product of $\|\hat{W}_\la^2(\rho) - \hat{W}_\la^2\|_{q/2}$
and a power of $\PP[ \hat{W}_\la \neq \hat{W}_\la(\rho) ] $. The first term is $O(\la^{2})$ by \eqref{bdd1} whereas the latter  is small by \eqref{bd2a}, the choice of $\rho$, and the arbitrariness of
$L$.

Likewise \eqref{diff2}  can be bounded by $\la^{-1}/2$ since
\begin{eqnarray*}
&&| \E [ \E  [\hat{W}_\la(\rho) |  \G_G]^2 - \E  [\hat{W}_\la |   \G_G]^2  ] |\\
&&=
|\E ( \E  [\hat{W}_\la(\rho) |   \G_G] +  \E  [\hat{W}_\la |   \G_G]  )( \E  [\hat{W}_\la(\rho) |   \G_G] -  \E  [\hat{W}_\la |   \G_G])|\\
&&
\leq C \la  \|  \E  [\hat{W}_\la(\rho) |  \G_G] - \E  [\hat{W}_\la |   \G_G] \|_2\\
&&\le C \la\sqrt{\E\left(\E\left((\hat{W}_\la(\rho)-\hat{W}_\la)^2\vert \G_G\right)\right)}=C \la\sqrt{\E(\hat{W}_\la(\rho)-\hat{W}_\la)^2},
\end{eqnarray*}
where the first inequality follows by the Cauchy-Schwarz inequality and Lemma~\ref{momentlemma}
and where the second inequality follows by the conditional Jensen inequality.  Using H\"older's inequality and the bound \eqref{bd2a},   we get that \eqref{diff2} is bounded by $\la^{-1}/2$, concluding the proof of \eqref{Xia14082901}.
The proofs of \eqref{Xia14082902} and part (b) follow the proof of (a) verbatim.   \qed

\vskip.3cm

\noindent{\em Proof of Theorem~\ref{Insurance}.}  We take $\tilde{S}_T:=[0,T]\times\D$ in Theorem~\ref{main-noti} and let $r$ be the length of $I$.
    {Let} $n(\rho,r,\tilde{S}_T)$ be the maximum number of {subsets $S_i \subset \tilde{S}_T$ of the form $(I + t_i) \times\D, t_i \in \R^+,$ in $\tilde{S}_T$ }  which are separated by $\eightrho$ with $\rho=\Theta(\ln T)$. 
     Then $\Vol_{{d + 1}}(\tilde{S}_T)=\Theta(T)$ and $n(\rho, {r}, \tilde{S}_T)=\Theta(T(\ln T)^{-1})$.
   Let  $\tilde{\P}_T^{\b \Psi}:= \P^{\b \Psi} \cap \tilde{S}_T$  in accordance with \eqref{tilP}.
       We show that \eqref{noti} is satisfied for all $S_i$, $1 \leq i \leq n(\rho, r, \tilde{S}_T)$ and then apply
     Theorem~\ref{main-noti} to  $\tilde{\P}_T^{\b \Psi}$.
Since the conditional distribution $ \tilde{W}_T|\tilde{\P}_T^{\b \Psi}\cap\{([0,T] \setminus I) \times \D \}$ is non-degenerate, we have
$$\E\Var[\tilde{W}_T|\tilde{\P}_T^{\b \Psi}\cap\{([0,T] \setminus I) \times \D \}]:= {d}_0>0. $$
 For $J\subset [0,T]\times\D$, we define
 $$
 M(J):=\int_{J}  \xi((t,\bs), \tilde{\P}_T^{\b \Psi} )
 {\bf 1} (D((t,\bs), \tilde{\H}_T^{\b \Psi}) \leq  \tworho) \tilde{\P}_T^{\b \Psi}(dt,d\bs).$$
 Then
  \bear
 &&\E\var [M(\tilde{S}_T)|\tilde{\P}_T^{\b \Psi}\cap\{\tilde{S}_T{\setminus} S_i\}]\nonumber\\
&=&\E\var[M(S_i^{2\rho})|\tilde{\P}_T^{\b \Psi}\cap \{\tilde{S}_T{\setminus} S_i\}]\nonumber\\
&=&\E\var[M((I\times \D)^{2\rho})|\tilde{\P}_T^{\b \Psi}\cap\{([0,T] \setminus I) \times \D \} ] \mbox{\hskip1cm(by translation invariance of $\xi$)}\nonumber\\
&=&\E\var[M(\tilde{S}_T)|\tilde{\P}_T^{\b \Psi}\cap \{([0,T] \setminus I) \times \D \}] \ge d_0-O(T^{-1}),\nonumber
\eear
where the inequality is due to Lemma~\ref{newvar}(b). Using Lemma~\ref{newvar}(b) again, we conclude that, for $T$ large,
$$\E\var[\tilde{W}_T|\tilde{\P}_T^{\b \Psi}\cap\{\tilde{S}_T{\setminus} S_i\}]\ge d_0-O(T^{-1}).$$
All conditions of Theorem~\ref{main-noti} are satisfied and it follows from \eqref{main1-02-noti} that
$$ d_K\left(  \frac{\tilde{W}_T - \E \tilde{W}_T}{\sqrt{\Var \tilde{W}_T} }  ,N(0,1)\right) =O((\ln T)^2 \Vol(\tilde{S}_T)n(\rho,r,\tilde{S}_T)^{-3/2})=O( (\ln T)^{3.5} T^{-1/2} ),$$
completing the proof. $\qed$

\begin{lemm} \label{varlb}  Assume that $\xi$ is translation invariant and the moment condition \eqref{mom}   holds for some $q \in (2, \infty)$.
Under conditions \eqref{ExpStab} and \eqref{assum2}
{there exist constants $\la_0\in (0, \infty)$ and $c\in (0, \infty)$ such that} 
for all $\la \in [\la_0, \infty)$ and all $\rho \in [c \ln \la, \infty)$ we have
\be \label{Lemlb} \Var [ W_\la(\rho)] \geq  c^{-1} b_0\la \rho^{-d}; \ \Var [ \hat{W}_\la(\rho)] \geq  c^{-1} b_0\la \rho^{-d}.\ee
\end{lemm}


\vskip.3cm

\noindent{\em Proof.} We only prove the first inequality as the second follows from identical methods.
Let $c\ge 2/c'$ such that Lemma~~\ref{newvar}(a) holds, where $c'$ is the constant such that  the cardinality of $I$ is $\lfloor c' \la/\rho^d \rfloor$.
 Let ${\cal F}$  be the sigma algebra generated by
${\H}^{\b \Psi} \cap  (\bigcup_{i\in I} Q_{i,r})^c$. By the conditional variance formula
$$
\Var [ W_\la(\rho)] = \Var[ \E [ W_\la(\rho) | \F] ] + \E \Var [W_\la(\rho) | \F]
\geq \E \Var [W_\la(\rho) | \F].
$$
Let $C_i:=\{x\in\R^d:\ d(x,Q_{i,r})\le \tworho\}$.
Then the $C_i$ are separated
by $\fourrho$ because the $Q_{i,r}$ are separated by at least $\eightrho$ (this is the reason why we chose the $\eightrho$ separation
in the first place). Also, the $C_i$ are contained in $Q_\la$.

For each $i \in I$ the sum
$\sum_{x \in {\H}_\la^{\b \Psi} \cap C_i } \xi(x,  {\H}_\la^{\b \Psi} { \setminusx}; \rho)$
depends on points distant at most 
$\tworho$ from $C_i$.
Thus the random variable $\E [W_\la(\rho)| \F]$ is a sum of
independent random variables since the $C_i$ are separated by $\fourrho$.
Thus we obtain

\bear\E \Var [W_\la(\rho) | \F] &=& \E \Var[\sum_{x \in {\H}_\la^{\b \Psi}} \xi(x,  {\H}_\la^{\b \Psi} { \setminusx}; \rho) | \F]\nonumber\\
&=& \E \sum_{i \in I} \Var [ \sum_{x \in {\H}_\la^{\b \Psi} \cap C_i} \xi(x,  {\H}_\la^{\b \Psi} { \setminusx}; \rho) | \F].\label{ana1}
\eear


Recall that $E^\epsilon = \{x \in \R^d: \ d(x, E) < \epsilon\}$ for any set $E$ and $\epsilon>0$.
For all $i \in I$,  the restrictions of $\F$ and $\F_i$ to $C_i^{\tworho}$ coincide.
For $x \in C_i$, we have that $\xi(x, \P_\la^{\beta \Psi}; \rho)$ depends only on points in $C_i^{\rho}$
and so we may thus replace
$\F$ with $\F_i$.   Since $\P_\la^{\beta \Psi}$ and $\P^{\beta \Psi}$ coincide on $C_i^{\rho}$ we may also replace $\xi(x, \P_\la^{\beta \Psi}; \rho)$ with
  $\xi(x, \P^{\beta \Psi}; \rho)$.   Also, we may replace the range of summation $x \in {\H}_\la^{\b \Psi} \cap C_i$ by
$x \in {\H}_\la^{\b \Psi}$ because the conditional sum
$$\sum_{x \in \P^{\b \Psi} \cap C_i^c \cap Q_\la} \xi(x,  {\H}^{\b \Psi} { \setminusx}; \rho)  | \F_i$$
is constant (indeed, if $x \in C_i^c$, then $ \xi(x,  {\H}^{\b \Psi} { \setminusx}; \rho)$ won't be affected by points in $Q_{i,r}$).

This yields
\be \label{ana2}
\E \Var [W_\la(\rho) | \F] =
\E \sum_{i \in I} \Var [ \sum_{x \in {\H}_\la^{\b \Psi} } \xi(x,  {\H}^{\b \Psi} { \setminusx}; \rho) | \F_i].
\ee


By Lemma~\ref{newvar}(a) for all $i \in I$,
$$ \E \Var [ \sum_{x \in {\H}_\la^{\b \Psi} } \xi(x, {\H}^{\b \Psi}{ \setminusx}; \rho){|\F}_i] \geq b_0/2.$$
 Thus
$$
\Var [ W_\la(\rho)] \geq \E \Var [W_\la(\rho) | \F] \geq \E \sum_{i \in I}  b_0/2 \geq  b_0 {c^{-1}}  \la \rho^{-d}.
$$
\qed

Roughly speaking, the factor $\la \rho^{-d}$ in \eqref{Lemlb}
is the cardinality of $I$, the index set of cubes of volume $r$, separated by $\eightrho$, and having
the property that the total score on each cube has positive variability.
For score functions which may not be translation invariant and/or are defined
on a subset $\tilde{S}_\la$ of $\R^d$, we have the following analog
of Lemma \ref{varlb}.  Recall the definition of $n(\rho, r,\tilde{S}_\la)$ right after \eqref{noti}.

\begin{lemm} \label{varlb-noti}  Assume  the moment condition \eqref{mom3}  holds for some $q \in (2, \infty)$.
Under conditions \eqref{ExpStab} and \eqref{noti}
{there exist constants $\la_0\in (0, \infty)$ and $c\in (0, \infty)$ such that}
for all $\la \in [\la_0, \infty)$ and all $\rho \in [c \ln \la, \infty)$ we have
 $$ \Var [ \tilde{W}_\la(\rho)] \geq  c^{-1} b_0 n(\rho, r, \tilde{S}_\la).$$
\end{lemm}

\noindent{\em Proof.} We follow the proof of {Lemma~\ref{varlb}}.
We write $\{\tilde{Q}_{i,r}: \ i\in \tilde{I}\}:={\cal C}(\rho,r, \tilde{S}_\la)$, the collection of cubes defined after \eqref{noti}.
 Let $\tilde{\F}_\la$  be the sigma algebra generated by
$\tilde{\H}_\la^{\b \Psi}\cap  (\bigcup_{i\in \tilde{I}} \tilde{Q}_{i,r})^c$. By the conditional variance formula
$$
\Var [ \tilde{W}_\la(\rho)] = \Var[ \E [\tilde{W}_\la(\rho) | \tilde{\F}_\la] ] + \E \Var [\tilde{W}_\la(\rho) | \tilde{\F}_\la]
\geq \E \Var [\tilde{W}_\la(\rho) | \tilde{\F}_\la].
$$
For $i \in \tilde{I}$, let $\tilde{C}_i:=\{x\in\tilde{S}_\la:\ d(x,\tilde{Q}_{i,r})\le \tworho\}$. 
Then the $\tilde{C}_i$ are separated
by $\fourrho$ because the $\tilde{Q}_{i,r}$ are separated by at least $\eightrho$.
 Also, the $\tilde{C}_i$ are contained in $\tilde{S}_\la$.

For each $i \in \tilde{I}$ the sum
$\sum_{x \in \tilde{{\H}}_\la^{\b \Psi} \cap \tilde{C}_i } \xi(x,  \tilde{{\H}}_\la^{\b \Psi} { \setminusx}; \rho)$
depends on points distant at most 
$\tworho$ from $\tilde{C}_i$.
Thus $\E [\tilde{W}_\la(\rho)| \tilde{\F}_\la]$ is a sum of
independent random variables since the $\tilde{C}_i$ are separated by $\fourrho$.
Thus we obtain the analog of \eqref{ana1}, namely
$$
\E \Var [\tilde{W}_\la(\rho) | \tilde{\F}_\la] =  \E \sum_{i \in \tilde{I}} \Var [ \sum_{x \in \tilde{{\H}}_\la^{\b \Psi} \cap \tilde{C}_i} \xi(x,  \tilde{{\H}}_\la^{\b \Psi} { \setminusx}; \rho) | \tilde{\F}_\la].
$$
Let $\tilde{\cal F}_{\la, i}$ be the sigma algebra generated by $\tilde{{\H}}_\la^{\b \Psi} \cap \tilde{Q}_{i,r}$.
For all $i \in \tilde{I}$,  the restrictions of $\tilde{\F}_\la$ and $\tilde{\F}_{\la, i}$ to $\tilde{C}_i^{\tworho}\cap \tilde{S}_\la$ coincide.

As in the proof of Lemma~\ref{varlb}, we obtain the analog of \eqref{ana2}, namely
$$
\E \Var [\tilde{W}_\la(\rho) | \tilde{\F}_\la] =
\E \sum_{i \in \tilde{I}} \Var [ \sum_{x \in \tilde{\H}_\la^{\b \Psi} } \xi(x,  \tilde{{\H}}^{\b \Psi} { \setminusx}; \rho) | \tilde{\F}_{\la,i}].
$$


If $\la \in [\la_0, \infty)$ and if
$\la_0$ is large enough, then by Lemma \ref{newvar}(b) for all $i \in \tilde{I}$,
$$ \E \Var [ \sum_{x \in \tilde{{\H}}_\la^{\b \Psi} } \xi(x, \tilde{{\H}}^{\b \Psi}{ \setminusx}; \rho){|\tilde{\F}}_{\la,i}] \geq b_0/2.$$
 Thus
$$
\Var [\tilde{W}_\la(\rho)] \geq \E \Var [\tilde{W}_\la(\rho) | \tilde{\F}_\la] \geq \E \sum_{i \in \tilde{I}}  b_0/2 \geq  b_0 \cdot \rm{card} (\tilde{I}).
$$
\qed

\vskip.1cm

\begin{lemm} \label{var630} If the moment condition \eqref{mom} holds for some $q \in (2, \infty)$ then $| \Var W_\la -\Var \hat{W}_\la| = o(\la).$
\end{lemm}

\noindent{\em Proof}. Put $\rho = c \ln \la$, $c$ large. By {\eqref{Xia14082902} and \eqref{Xia14082901} with $G=\emptyset$ we have $| \Var W_\la(\rho) -\Var W_\la| = o(1)$
and $| \Var \hat{W}_\la(\rho) -\Var \hat{W}_\la| = o(1)$.}  So it is enough to prove
$| \Var W_\la(\rho) -\Var \hat{W}_\la(\rho)| = o(\la).$  We have
$$| \Var W_\la(\rho) - \Var \hat{W}_\la(\rho)| \leq  \Var( W_\la(\rho) - \hat{W}_\la(\rho)) + 2 {\rm cov}( W_\la(\rho) - \hat{W}_\la(\rho),
\hat{W}_\la(\rho)).$$
The scores $\xi(x,  \P^{\b \Psi}_\la; \rho)$ and $\xi(x,  \P^{\b \Psi}; \rho)$ coincide  when $x \in Q_\la$ is distant at least
$\rho$ from $\partial Q_\la$. Thus  $W_\la(\rho) - \hat{W}_\la(\rho) = U_\la - V_\la$, where
$$
U_\la := \sum_{x \in  \P^{\b \Psi}_\la \cap (\partial Q_\la)^\rho } \xi(x,  \P^{\b \Psi}_\la; \rho); \ \ V_\la
:= \sum_{x \in  \P^{\b \Psi} \cap (\partial Q_\la)^\rho } \xi(x,  \P^{\b \Psi}; \rho).$$
 Lemma~\ref{momentlemma} with $ q = 2$ and $q' > 2$ ensures
$\Var U_\la$ and $\Var V_\la$ are both of {order $O( (\Vol (\partial Q_\la)^\rho)^2 )$.}
These bounds and the  formula $\Var[U_\la - V_\la] = \Var U_\la + \Var V_\la - 2 \Cov[U_\la, V_\la]$ shows that
$\Var[U_\la - V_\la] = o(\la)$.
By the Cauchy-Schwarz inequality and Lemma \ref{momentlemma}, we obtain ${\rm cov}[ W_\la(\rho) - \hat{W}_\la(\rho),
\hat{W}_\la(\rho)] = o(\la)$ as well.   \qed

\vskip.5cm
 We need one more lemma.  It shows that if fluctuations  of $\hat{W}_\la$  are not of volume order then they are necessarily at most of surface order and vice versa.  A version of this dichotomy appears in the statistical physics literature \cite{MY} and
 also in \cite{BDY}.    We do not have any natural examples of $\hat{W}_\la$ which are defined on all of $Q_\la$ and which have fluctuations  at
 most of surface order.
 However, when ancestor clans and stabilization radii have slowly decaying tails we expect that $\Var \hat{W}_\la$ behaves less like
 a sum of i.i.d. random variables and more like a sum of random variables with very long range dependencies, presumably giving rise to
 smaller 
 fluctuations. When the score at $x$ is  allowed to depend on nearby point configurations as well as on nearby scores, then Martin and Yalcin
 \cite{MY} establish conditions giving surface order fluctuations.

\begin{lemm}
\label{varld}  Let $\xi$ be translation invariant.  Either $\Var \hat{W}_\la = \Omega(\la)$ or
$\Var \hat{W}_\la = O(\la^{(d-1)/d}).$
\end{lemm}

\noindent{\em Proof}. 
Recall the definitions of $c^\xi(x)$ and $c^\xi(x,y)$ at \eqref{defcx} and \eqref{defc2x}, respectively.
Similar to the proof of Theorem 2.2 of \cite{SY}, by the integral characterization of Gibbs point processes, as in Chapter 6.4 of \cite{MW},
it follows from the Georgii-Nguyen-Zessin formula that
$$
\Var \hat{W}_\la = \Var \sum_{x \in \P^{\b \Psi}_\la} \xi(x, {\H}^{\b \Psi} \setminusx) = {\tau} \int_{Q_\la} c^{\xi^2}(x) dx - {\tau}^2\int_{Q_\la} \int_{Q_\la}c^\xi(x,y) dy dx.
$$
Note that $c^\xi(x,y)$ decays exponentially fast with $|x - y|$, as shown in Lemmas 3.4 and 3.5 of
\cite{SY}. By translation invariance of $\xi$ and stationarity of $\P^{\b \Psi}$ we get
\be \label{E1}
\Var \hat{W}_\la = {\tau}c^{\xi^2}(\0)\la - {\tau}^2 \int_{Q_\la} \int_{\R^d}c^\xi(\0,y-x) {\bf 1}( y \in Q_\la)  dy dx
\ee
$$
= {\tau}c^{\xi^2}(\0)\la- {\tau}^2 \int_{Q_\la} \int_{\R^d}c^\xi(\0,y) {\bf 1}( x + y \in Q_\la)  dy dx:= I_\la + II_\la.
$$
Now
$$
\la^{-1} II_\la  = -{\tau}^2 \la^{-1} \int_{Q_\la} \int_{\R^d}c^\xi(\0,y) {\bf 1}( x \in Q_\la - y)  dy dx$$
and writing ${\bf 1}( x \in Q_\la - y)$ as $1 -  {\bf 1}( x \in (Q_\la - y)^c)$ gives
$$
\la^{-1} II_\la = -{\tau}^2 \int_{\R^d} c^\xi(\0,y) dy + \la^{-1} {\tau}^2\int_{\R^d} \int_{Q_\la} c^\xi(\0,y) {\bf 1}( x \in \R^d \setminus (Q_\la - y))  dxdy.$$
As in \cite{MY}, for all $y \in \R^d$, put $\gamma_{Q_\la}(y):= \Vol_d(Q_\la \cap (\R^d \setminus (Q_\la - y))).$  Then
\be \label{E2}
\la^{-1} \Var \hat{W}_\la =  \la^{-1} I_\la + \la^{-1} II_\la = {\tau} c^{\xi^2}(\0) - {\tau}^2 \int_{\R^d} c^\xi(\0,y) dy +  \la^{-1} {\tau}^2 \int_{\R^d} c^\xi(\0,y)  \gamma_{Q_\la}(y)  dy.\ee
Now we assert that
\be \label{E3}
\lim_{\la \to \infty} \la^{-1} \int_{\R^d} c^\xi(\0,y)   \gamma_{Q_\la}(y)  dy = 0.\ee
Indeed, {by Lemma 1 of \cite{MY},} we have  $ \la^{-1}  \gamma_{Q_\la}(y) \to 0$ and since
 $ \la^{-1} c^\xi(\0,y)  \gamma_{Q_\la}(y) $ is dominated by $c^\xi(\0,y)$, which decays exponentially fast,
 the result follows by the dominated convergence theorem.  
Collecting terms in \eqref{E1}-\eqref{E3} and recalling \eqref{defsig} gives
\be \label{hatlimit}
\lim_{\la \to \infty} \la^{-1} \Var \hat{W}_\la = {\tau} c^{\xi^2}(\0) - {\tau}^2\int_{\R^d} c^\xi(\0,y) dy = {\tau} \sigma^2(\xi, {\tau}) \in [0, \infty),
\ee
where we note $\sigma^2(\xi,\tau)$ is finite by the exponential decay of $c^\xi(\0,y)$
as shown in Lemma 3.5 of \cite{SY}.

It follows that if $\Var \hat{W}_\la$ is not
of volume order then we have ${\tau} c^{\xi^2}(\0) - {\tau}^2\int_{\R^d} c^\xi(\0,y) dy = 0$.  Using this identity in
\eqref{E2}, multiplying \eqref{E2} by $\la^{1/d}$, and taking limits gives
\be \label{E4}
\lim_{\la \to \infty} \la^{-(d-1)/d} \Var \hat{W}_\la = \lim_{\la \to \infty} {\tau}^2  \la^{-(d-1)/d} \int_{\R^d} c^\xi(\0,y)   \gamma_{Q_\la}(y)  dy. \ee
Now as in \cite{MY}, we have $\la^{-(d-1)/d} \gamma_{Q_\la}(y) \leq C |y|,$
showing that the integrand in \eqref{E4} is dominated by an integrable function.
By Lemma 1 of \cite{MY}, there is a function $\gamma: \R^d \to \R^+$ such that
$$
\lim_{\la \to \infty} \la^{-(d-1)/d} \gamma_{Q_\la}(y) = \gamma(y).
$$
By dominated convergence we get the desired result:
$$
\lim_{\la \to \infty} \la^{-(d-1)/d} \Var \hat{W}_\la =  {\tau}^2 \int_{\R^d} c^\xi(\0,y) \gamma(y) dy < \infty,
$$
where once again the integral is finite by the exponential decay of $c^\xi(\0,y)$.
\qed

\section{Proofs of Theorems \ref{main0}-\ref{main-noti}} \label{proofofmainthm}
\allco


\noindent{\em Proof of Theorem~\ref{main0}.}  Combining  \eqref{hatlimit} and  Lemma \ref{var630} we obtain $\liml \la^{-1}  \Var {W}_\la =
{\tau}\sigma^2(\xi, {\tau}),$
giving \eqref{Var11}.  Now assume non-degeneracy \eqref{assum2} and put $\rho = {c} \ln \la$.  By Lemma
\ref{varlb} we have
$$
\lim_{\la \to \infty} \la^{-(d-1)/d} \Var \hat{W}_\la(\rho) = \infty
$$
and therefore by {\eqref{Xia14082901} with $G=\emptyset$} we have
$\lim_{\la \to \infty} \la^{-(d-1)/d}\Var \hat{W}_\la = \infty.
$
By Lemma \ref{varld} we have $\Var \hat{W}_\la = \Omega(\la)$ and Lemma \ref{var630}
gives $\sigma^2(\xi, {\tau}) > 0$, as desired. \qed

\vskip.3cm

\noindent{\em Proof of Theorem~\ref{main}.} We use a result based on the Stein method to derive rates of normal convergence. We  follow the set-up of \cite{BX1}, as this yields rates which are a slight improvement
over the methods of \cite{SY}.  
Given an admissible  Gibbs point process ${\H}_\la^{\b \Psi}$ with both $\beta$ and
$\Psi$ fixed, we shall simply write ${\H}_\la$ for ${\H}_\la^{\b \Psi}$.
 Our first goal
is to get rates of normal convergence for $W_\la(\rho)$ defined at \eqref{Wla}. Then we use this to obtain rates
for $W_\la$.  Without loss of generality, we assume $p \in (2, q)$    
  and we show for all $\rho \in (0, \infty)$: 
\begin{equation}d_K\left(\frac{W_\la(\rho)-\E W_\la(\rho)}{\sqrt{\var(W_\la(\rho))}},N(0,1)\right)
=O\left( ( \var W_\la(\rho))^{-p/2}\la w_{q}^p{\rho^{d(p-1)}}+ (\var W_\la(\rho))^{-1/2}w_{q}\rho^d\right)\label{proofmainthm00-1}\end{equation}
and, if \eqref{assum2} holds and if  \eqref{mom} holds for some $q \in (3, \infty)$,
\begin{equation}d_K\left(\frac{W_\la(\rho)-\E W_\la(\rho)}{\sqrt{\var(W_\la(\rho))}},N(0,1)\right)
=O\left( \rho^{{2d}} \la^{-1/2} \right).\label{proofmainthm00-2}\end{equation}

The proof goes as follows. 
 The local dependence
condition LD3 of \cite{BX1} requires
 for each $x\in Q_\la$  three nested neighborhoods $A_x$, $B_x$ and $C_x$ which satisfy $B_r(x)\subset A_x\subset B_x\subset C_x$ as $r\downarrow 0$
  and such that the sum of scores over points in $B_r(x)$ (resp. $A_x$, $B_x$) are independent of the sum of scores over points in $(A_x^r)^c$ (resp. $B_x^c$, $C_x^c$). We claim that  $W_\la(\rho)$ satisfies the local dependence
condition LD3 with the neighborhoods $A_x:=B_{2\rho}(x)$,
$B_x:=B_{4\rho}(x)$ and $C_x:=B_{6\rho}(x)$,
$x\in Q_\la$. Indeed, this follows immediately since $\xi(\cdot,  {\H}_\la^{\b \Psi} \setminus {\{\cdot\}};
\rho)$
 enjoys spatial independence over
sets separated by more than $2\rho$, as already noted in the discussion after \eqref{Wla}.

It follows from Corollary 2.2 of \cite{BX1} that
$$d_K\left(\frac{W_\la(\rho)-\E W_\la(\rho)}{\sqrt{\var(W_\la(\rho))}},N(0,1)\right)\le 48\varepsilon_3+160\varepsilon_4+2\varepsilon_5,$$
where, with $R(dx):= | \xi(x,  {\H}_\la; \rho)| {\H}_\la(dx)$, $N(C_x):= B_{10\rho}(x)$, and $p \in (2,
\infty)$,
 \bean
\varepsilon_3&:=& (\var W_\la(\rho))^{-p/2}\E\int_{Q_\la} R(N(C_x))^{p-1}R(dx),\\
\varepsilon_4&:=& (\var W_\la(\rho))^{-p/2}\int_{Q_\la} \E R(N(C_x))^{p-1}\E R(dx),\\
\varepsilon_5&:=& (\var W_\la(\rho))^{-1/2}\sup_{x \in Q_\la} \E R(N(C_x)).\\
\eean
{We write} $G_{x,\la} := \{ {D(x,{\H}_\la) \leq \rho}\}$.
For $\varepsilon_3$, we have by definition of $R(dx)$ 
 that \bear 
 &&\E\int_{Q_\la}  R(N(C_x))^{p-1}R(dx) 
 \nonumber\\
 & = &  \E\int_{Q_\la}\left(\int_{N(C_x)}|\xi(z,  {\H}_\la \setminusz )|{\bf 1}(G_{z,\la}){\H}_\la(dz)\right)^{p-1} |\xi(x,
 {\H}_\la \setminusx)|{\bf 1}(G_{x,\la}){\H}_\la(dx) \nonumber\\
 & \le &  \E\int_{Q_\la}\left(\int_{N(C_x)}|\xi(z,  {\H}_\la \setminusz )|{\H}_\la(dz)\right)^{p-1} |\xi(x,
 {\H}_\la \setminusx)|{\H}_\la(dx). \nonumber\eear

\noindent H\"older's inequality $(\int_D |f| \mu(dx))^{p-1} \leq \int_D
|f|^{p-1} \mu(dx) \cdot \mu(D)^{p-2}$ gives  that \bear
&& \E\int_{Q_\la}  R(N(C_x))^{p-1}R(dx)  \nonumber\\
&\le&
 \E\int_{Q_\la}\int_{N(C_x)}|\xi(z, {\H}_\la \setminusz )|^{p-1}{\H}_\la(dz) \cdot {\H}_\la(N(C_x))^{p-2}
|\xi(x,  {\H}_\la\setminusx)|{\H}_\la(dx)\nonumber\\
&\le& \E\int_{Q_\la}|\xi(x,  {\H}_\la \setminusx)|^{p}{\H}_\la(N(C_x))^{p-2}
{\H}_\la(dx)\nonumber\\
&&+  \E\int_{Q_\la}\int_{N(C_x){\setminusx}}|\xi(z,  {\H}_\la \setminusz)|^{p-1}{\H}_\la(dz){\H}_\la (N(C_x))^{p-2}
|\xi(x,  {\H}_\la\setminusx)|{\H}_\la(dx) \nonumber,
\eear
where we write $\int_{N(C_x)} \cdots {\H}_\la(dz)$ as $\int_{\{x\}} \cdots {\H}_\la(dz) + \int_{N(C_x) \setminus \{x\}}\cdots {\H}_\la(dz)$. The inequality $|a\|b|^{p-1}\le |a|^p+|b|^p$
gives
\bear
&& \E\int_{Q_\la}  R(N(C_x))^{p-1}R(dx)  \nonumber\\
&\le& \E\int_{Q_\la}|\xi(x,  {\H}_\la \setminusx)|^{p}{\H}_\la(N(C_x))^{p-2}{\H}_\la(dx)\nonumber\\
&&+ \E\int_{Q_\la}\int_{N(C_x){\setminusx}}\left(|\xi(z, {\H}_\la\setminusz)|^{p}+|\xi(x,  {\H}_\la\setminusx)|^{p}\right) \cdot{\H}_\la(N(C_x))^{p-2} {\H}_\la(dz){\H}_\la(dx).\nonumber\eear
Splitting the last integral into two
integrals gives \bear
&& \E\int_{Q_\la}  R(N(C_x))^{p-1}R(dx)  \nonumber\\
 &\le&
\E\int_{Q_\la}|\xi(x,  {\H}_\la  \setminusx)|^{p}{\H}_\la(N(C_x))^{p-2}{\H}_\la(dx)\nonumber\\
&&+ \E\int_{Q_\la}\int_{N(C_x){\setminusx}}|\xi(z,  {\H}_\la  \setminusz)|^{p}{\H}_\la(N(C_x))^{p-2}
{\H}_\la(dz){\H}_\la(dx)\nonumber\\
&&+ \E\int_{Q_\la}|\xi(x,  {\H}_\la  \setminusx)|^{p}{\H}_\la(N(C_x))^{p-1}{\H}_\la(dx)\nonumber\\
&\le& \E\int_{Q_\la}{\H}_\la(N(C_x))^{p-2}|\xi(x,  {\H}_\la  \setminusx)|^{p}{\H}_\la(dx)\nonumber\\
&&+ \E\iint_{0<d(x,z)\le 10 \rho}|\xi(z,  {\H}_\la  \setminusz)|^{p}{\H}_\la(N(C_x))^{p-2}
{\H}_\la(dx){\H}_\la(dz)\nonumber\\
&&+ \E\int_{Q_\la}{\H}_\la(N(C_x))^{p-1}|\xi(x,  {\H}_\la  \setminusx)|^{p}{\H}_\la(dx).\nonumber\eear Now
integrating the double integral gives
\bear
&& \E\int_{Q_\la}  R(N(C_x))^{p-1}R(dx)  \nonumber\\
&\le&
\E\int_{Q_\la}{\H}_\la(N(C_x))^{p-2}|\xi(x,  {\H}_\la  \setminusx)|^{p}{\H}_\la(dx)\nonumber\\
&&+ \E\int_{Q_\la}|\xi(z,  {\H}_\la  \setminusz)|^{p} \cdot {\H}_\la(B_{20\rho}(z))^{p-1}{\H}_\la(dz)\nonumber\\
&&+\E\int_{Q_\la}{\H}_\la(N(C_x))^{p-1}|\xi(x,  {\H}_\la\setminusx)|^{p}{\H}_\la(dx).\nonumber
\eear
Combining integrals and using H\"older's inequality for $p_1 \in (1, q/p)$ gives \bear
&& \E\int_{Q_\la}  R(N(C_x))^{p-1}R(dx)  \nonumber\\
&\le&3\E\int_{Q_\la}|\xi(z,  {\H}_\la  \setminusz)|^{p}{\H}_\la(B_{20\rho}(z))^{p-1}{\H}_\la(dz)\nonumber\\
&\le&3\left\{\E\int_{Q_\la} {\H}_\la(B_{20\rho}(z))^{\frac{(p-1)p_1}{p_1-1}}{\H}_\la(dz)\right\}^{\frac{p_1-1}{p_1}}\left\{\E\int_{Q_\la} |\xi(z,  {\H}_\la \setminusz)|^{pp_1}{\H}_\la(dz)\right\}^{\frac1{p_1}}.\label{proofmainthm01} \eear
{Since ${\H}_\la^{\b \Psi}$ is a Gibbs point process, we apply the Georgii-Nguyen-Zessin
integral characterization of Gibbs point processes {\cite{MW}} to see that the
conditional probability of observing an extra point of ${\H}_\la^{\b \Psi}$ in the volume element $dz$, given that configuration without
that point, equals
 $\exp( - \beta \Delta^{{\Psi}}( \{z\}, {\H}_\la^{\b \Psi}
) )dz\le dz$, where $\Delta^{{\Psi}}( \{z\}, {\H}_\la^{\b \Psi})$ is defined at \eqref{Dell}.
Using that ${\E}\H_\la^{\b \Psi}(dx) \leq \tau dx$,}
we have from \eq{proofmainthm01} that
\bear
&&\E\int_{Q_\la}  R(N(C_x))^{p-1}R(dx)\nonumber\\
&\le&3 \tau\left\{\E\int_{Q_\la} \left({\H}_\la(B_{20\rho}(z))+1\right)^{\frac{(p-1)p_1}{p_1-1}} dz\right\}^{\frac{p_1-1}{p_1}}\left\{\E\int_{Q_\la} |\xi({z},  {\H}_\la \cup \{z\})|^{pp_1}
dz\right\}^{\frac1{p_1}}.\label{proofmainthm02} \eear
Notice that ${\H}_\la(B_{20\rho}(x))$ is stochastically bounded by ${{\rm
Po}(\tau M)}$
with $M:= \mbox{Vol}(B_{20\rho}(0))$, we have
from Lemma~4.3 of \cite{BX1} that $\E\{{\H}_\la(B_{20\rho}(x))+1\}^{(p-1)p_1/(p_1-1)}\le
{c_1}{\rho^{d(p-1)p_1/(p_1-1)}}$, giving
$$ \varepsilon_3\le3 \tau \var(W_\la(\rho))^{-p/2}
{c_1}^{\frac{p_1-1}{p_1}}{\rho^{d(p-1)}}\la^{(p_1-1)/p_1}\left\{\E\int_{Q_\la} |\xi(x,  {\H}_\la \cup \{x\})|^{pp_1} dx\right\}^{\frac1{p_1}}.
$$
Then since  $w_{p p_1} \leq w_q$, we have \beqn \varepsilon_3\le3\tau \la\var(W_\la(\rho))^{-p/2}
{c_1}^{\frac{p_1-1}{p_1}}w_{q}^p{\rho^{d(p-1)}}.\label{proofmainthm03} \eeqn

Next, we bound $\varepsilon_4$. To this end, let
$p_2:= pp_1/(p-1)$, we again replace the indicator function with $1$ and then apply H\"older's inequality to get
\bear
&&\int_{Q_\la} \E R(N(C_x))^{p-1}\E R(dx)
\nonumber\\
&=& \int_{Q_\la}\E\left(\int_{N(C_x)}|\xi(z,  {\H}_\la  \setminusz)|{\bf
1}(G_{z,\la}){\H}_\la(dz)\right)^{p-1}
\E|\xi(x,  {\H}_\la  \setminusx)|{\bf 1}(G_{x,\la}){\H}_\la(dx)\nonumber\\
&{\le}& \int_{Q_\la}\E\left(\int_{N(C_x)}|\xi(z,  {\H}_\la  \setminusz)|{\H}_\la(dz)\right)^{p-1}
\E|\xi(x,  {\H}_\la  \setminusx)|{\H}_\la(dx)\nonumber\\
&\le& \int_{Q_\la}\E\left\{\int_{N(C_x)}|\xi(z,  {\H}_\la \setminusz)|^{p-1}{\H}_\la(dz){\H}_\la(N(C_x))^{p-2}\right\}
\E|\xi(x,  {\H}_\la  \setminusx)|{\H}_\la(dx)\nonumber\\
&\le& \int_{Q_\la}\E\left\{\int_{N(C_x)}|\xi(z,  {\H}_\la  \setminusz)|^{p-1}{\H}_\la(B_{20\rho}(z))^{p-2}{\H}_\la(dz)\right\} \E|\xi(x,  {\H}_\la  \setminusx)|{\H}_\la(dx)\nonumber\\
&\le&\int_{Q_\la}\left\{\E\int_{N(C_x)}|\xi(z,  {\H}_\la \setminusz)|^{p_2(p-1)}{\H}_\la(dz)\right\}^{\frac1{p_2}}\nonumber\\
&&\left\{\E\int_{N(C_x)}{\H}_\la(B_{20\rho}(z))^{(p-2)\frac{p_2}{p_2-1}}{\H}_\la(dz)\right\}^{\frac{p_2-1}{p_2}}
\E|\xi(x,  {\H}_\la  \setminusx)|{\H}_\la(dx).\label{proofmainthm04}
\eear Reasoning as for \eq{proofmainthm02}, we obtain from
\eq{proofmainthm04} that \bean
&&\int_{Q_\la} \E R(N(C_x))^{p-1}\E R(dx)\nonumber\\
&\le& \int_{Q_\la}\left\{\int_{N(C_x)}\E|\xi(z,  {\H}_\la \cup \{z\})|^{p_2(p-1)} \tau dz\right\}^{\frac1{p_2}}\nonumber\\
&&\left\{\int_{N(C_x)}\E\left({\H}_\la(B_{20\rho}(z))+1\right)^{(p-2)\frac{p_2}{p_2-1}}\tau dz\right\}^{\frac{p_2-1}{p_2}}
\E|\xi(x,  {\H}_\la  \setminusx)|{\H}_\la(dx)\nonumber\\
&\le&  \tau^2 w_{pp_1}^{p-1}{c_2}^{\frac{p_2-1}{p_2}}{\rho^{d(p-2)}}\int_{Q_\la}\left\{\int_{N(C_x)}
dz\right\}^{\frac1{p_2}} \left\{\int_{N(C_x)}dz\right\}^{\frac{p_2-1}{p_2}}
w_{pp_1} dx\nonumber\\
&\le&w_{pp_1}^{p}{c_3}\la{\rho^{d(p-1)}}. \eean Hence \beqn\varepsilon_4\le
(\var W_\la(\rho))^{-p/2}w_{q}^{p}{c_3}\la{\rho^{d(p-1)}},\label{proofmainthm05}
\eeqn
showing that the bounds for $\varepsilon_3$ and $\varepsilon_4$ coincide. Turning to  $\varepsilon_5$, we have
\bear
\varepsilon_5 &\le & (\var W_\la(\rho))^{-1/2}\sup_{x \in Q_\la} \E\left(\int_{N(C_x)}|\xi(z,  {\H}_\la \setminusz)|{\H}_\la(dz)\right)\nonumber\\
&\le &(\var W_\la(\rho))^{-1/2}\sup_{x \in Q_\la} \left(\int_{N(C_x)}\E|\xi(z,  {\H}_\la \cup \{z\})|\tau dz\right)\nonumber\\
&\le&\var(W_\la(\rho))^{-1/2}\sup_{x \in Q_\la} \left(\int_{N(C_x)}\left\{\E|\xi(z,  {\H}_\la \cup\{z\})|^{pp_1}\right\}^{\frac1{pp_1}} \tau dz\right)\nonumber\\
&\le&\var(W_\la(\rho))^{-1/2}w_{q}{c_4}\rho^d.\label{proofmainthm06} \eear
Combining estimates \eq{proofmainthm03}, \eq{proofmainthm05} and
\eq{proofmainthm06}, we get \Ref{proofmainthm00-1}.

Assuming condition \eqref{mom}, using \eqref{Xia14082902} with $G=\emptyset$ and Theorem~\ref{main0}, we have
$\var[W_\la(\rho)] \geq {c_5}  \la.$
When $p = 3$, {this, together with \Ref{proofmainthm00-1},} gives {\Ref{proofmainthm00-2}.}

To complete the proof, we need to replace $W_\la(\rho)$ with $W_\la$. We rely heavily on {Lemma~\ref{newvar}} for this.
Note for all $\epsilon_1\in\R$ and $\epsilon_2>-0.6$,
\bear
&&d_K(N(0,1),N(\epsilon_1,1+\epsilon_2))\le d_K(N(0,1),N(\epsilon_1,1))+d_K(N(\epsilon_1,1),N(\epsilon_1,1+\epsilon_2))\nonumber\\
&&\le\frac{|\epsilon_1|}{\sqrt{2\pi}}+\frac{|\epsilon_2|}{\sqrt{2e\pi}}.\label{normal} \eear

Now $d_K(X, N(0,1)) = d_K(aX, N(0, a^2)) = d_K(aX + b, N(b, a^2))$ holds for $X$ with $\E X=0$ and all constants $a$ and $b$.  Hence
\bear
&&d_K\left(\frac{W_\la-\E W_\la}{\sqrt{\var W_\la}},N(0,1)\right)=d_K\left(\frac{W_\la-\E W_\la(\rho)}{\sqrt{\var W_\la(\rho)}},N\left(\frac{\E W_\la-\E W_\la(\rho)}{\sqrt{\var W_\la(\rho)}},\frac{\var W_\la}{\var W_\la(\rho)}\right)\right)\nonumber\\
&\le&d_K\left(\frac{W_\la-\E W_\la(\rho)}{\sqrt{\var W_\la(\rho)}},N(0,1)\right)+d_K\left(N(0,1),N\left(\frac{\E W_\la-\E W_\la(\rho)}{\sqrt{\var W_\la(\rho)}},\frac{\var W_\la}{\var W_\la(\rho)}\right)\right)\nonumber\\
\eear
by the triangle inequality for $d_K$. 
{ Now for any random variables $Y$ and $Y'$ we have
\be \label{ax}
d_K(Y, N(0,1)) \leq d_K(Y', N(0,1)) + \PP[Y \neq Y']
\ee
which follows from $| \PP[Y \leq t] - \Phi(t)| \leq | \PP[Y' \leq t] - \Phi(t)| + |\PP[Y' \leq t] - \PP[Y \leq t]|$.   We have by \eqref{ax}
and \eqref{normal} that}
\bear
&&  d_K\left(\frac{W_\la-\E W_\la}{\sqrt{\var W_\la}},N(0,1)\right) \nonumber\\
&\le&\PP[W_\la\ne W_\la(\rho)]+d_K\left(\frac{W_\la(\rho)-\E W_\la(\rho)}{\sqrt{\var W_\la(\rho)}},N(0,1)\right)\nonumber\\
&&+\frac{1}{\sqrt{2\pi}}\left|\frac{\E W_\la-\E W_\la(\rho)}{\sqrt{\var W_\la(\rho)}}\right|+\frac{1}{\sqrt{2e\pi}}\left|\frac{\var W_\la-\var W_\la(\rho)}{\var W_\la(\rho)}\right|.\label{Xia14082903}
\eear
However,
the Cauchy-Schwarz inequality ensures
$$|\E W_\la-\E W_\la(\rho)|\le \|W_\la-W_\la(\rho)\|_2\PP(W_\la\ne W_\la(\rho))^{1/2}\le \la^{-1},$$
where the last inequality is due to \eqref{bdd1}, \eqref{bd2b} and the arbitrariness of
$L$. Hence, it follows from \eqref{Xia14082903} that
$$  d_K\left(\frac{W_\la-\E W_\la}{\sqrt{\var W_\la}},N(0,1)\right)
\le\la^{{-2}} +{O\left( ( \var W_\la)^{-p/2}\la(\ln\la)^{d(p-1)}\right)},
$$
where we use \eqref{bd2b} with $L=2$, \eqref{proofmainthm00-1} and {\eqref{Xia14082902} with $G=\emptyset$}. \qed

\vskip.5cm

\noindent{\em Proof of  Theorem \ref{main-noti}}. The bound \eqref{VartW} follows from Lemma~\ref{varlb-noti} and Lemma~\ref{newvar}(b) with $G=\emptyset$.
{ The proof of \eqref{main1-01-noti} 
follows by replacing $Q_\la$ with
 $\tilde{S}_\la$ in the proof of \eqref{main1-01}, whereas  \eqref{main1-02-noti} follows by combining \eqref{main1-01-noti}
 and \eqref{VartW}.} \qed


\vskip.5cm

\noindent{\bf Acknowledgements.} We thank Yogeshwaran Dhandapani for
showing us Lemma \ref{varld}, {which essentially first appeared in \cite{BDY}}.  J. Yukich gratefully acknowledges generous and kind support from the Department of Mathematics and Statistics
at the University of Melbourne, where this work was initiated.

Aihua Xia, Department of Mathematics and Statistics, The University of Melbourne, Parkville, VIC 3010: {\texttt aihuaxia@unimelb.edu.au}

J. E. Yukich, Department of Mathematics, Lehigh University,
Bethlehem PA 18015:
\\
{\texttt joseph.yukich@lehigh.edu}

\end{document}